\documentclass[11pt,reqno]{amsart}
\usepackage{amsfonts, amsthm, amssymb, amsmath}
\usepackage{cite}

\theoremstyle{plain}
\newtheorem{thm}{Theorem}[section]
\newtheorem{lem}{Lemma}[section]
\newtheorem{cor}{Corollary}[section]
\newtheorem{claim}{Claim}[section]

\theoremstyle{definition}
\newtheorem{defn}{Definition}[section]
\newtheorem{rem}{Remark}[section]
\newtheorem{prob}{Problem}

\def\restrict#1{|_{#1}}
\def\ind{\operatorname{ind}}
\newcommand{\ord}{\mathop{\mathrm{ord}}}
\newcommand{\len}{\mathop{\mathrm{len}}}
\newcommand\aaa{{\boldsymbol{a}}}
\newcommand\bbb{{\boldsymbol{b}}}
\newcommand\ccc{{\boldsymbol{c}}}
\newcommand\ddd{{\boldsymbol{d}}}
\newcommand\xxx{{\boldsymbol{x}}}
\newcommand\yyy{{\boldsymbol{y}}}
\newcommand\zzz{{\boldsymbol{z}}}
\def\J_#1(#2,#3){J_{#1}^{#3}(#2)}
\def\circledsum{\rlap{\kern6.5pt$\circ$}\sum}
\def\dcircledsum{\rlap{\kern8.8pt$\circ$}\sum}

\begin{document}

\title{Nonequality of Dimensions for Metric Groups}

\author{Ol'ga~V.~Sipacheva}

\email{osipa@rol.ru}

\address{Department of General Topology and Geometry, 
Mechanics and Mathematics Faculty, Moscow State University, 
Moscow 119992, Russia}

\thanks{This work was financially supported by the Russian Foundation 
for Basic Research (project no.~06-01-00761).}

\begin{abstract}
An embeddability criterion
for zero-dimensional 
metrizable topological spaces in zero-dimensional 
metrizable topological groups is given.  
A space which can be embedded as a closed subspace in a 
zero-dimensional metrizable group 
but is not strongly zero-dimensional is 
constructed; thereby, an example of a metrizable group with 
noncoinciding dimensions $\ind$ and $\dim$ is obtained. It is 
proved that one of Kulesza's 
zero-dimensional metrizable spaces cannot be embedded in 
a metrizable zero-dimensional group.  
\end{abstract}



\maketitle

The presence of a topological group structure on a topological space 
has a strong influence on many properties of the space;  
a classical illustration is the metrizability of 
any first countable topological group. 
The dimensional properties are no exception. 
Thus, $\ind G= \dim G =\operatorname{{Ind}}G$ for any locally 
compact group $G$~\cite{Pas} and $\ind G= \operatorname{{Ind}}G$ 
for any topological group $G$ which is a  Lindel\"of
$\Sigma$-space~\cite{Shakh}, while for a general 
topological space, these three dimensions can be pairwise different,
even if the space is compact~\cite{Fil}. 

The purpose of this paper is to investigate the dimensional properties 
of metrizable topological groups. The celebrated theorem of 
Kat\v etov~\cite{Kat} says 
that $\dim X=\operatorname{{Ind}}X$ for any metric 
space $X$; however, there exist examples 
of metrizable spaces with noncoinciding dimensions $\ind$ and 
$\dim$. The first (very involved) example of such a space 
was constructed by Roy in 1968~\cite{Roy}. Since then, much simpler 
examples with various 
additional properties have been suggested (see, e.g., 
\cite{Kul, Kul2, Kul3, Mro1, 1, Ost}), 
but the question about the coincidence of dimensions for 
metrizable topological groups has remained open (apparently, for the 
first time, it was stated by Mishchenko in 1964~\cite{Mishch}). 

In the first section of this paper, we prove a criterion
for the embeddability of zero-dimensional metrizable topological 
spaces in zero-dimensional metrizable topological groups. 
This criterion was formulated by Mishchenko in~\cite{Mishch}, but its 
proof has never been published; Mishchenko himself confessed to this 
author in a private communication that he had retained neither 
notes nor recollections of the proof.  The spaces embeddable in 
zero-dimensional topological groups occupy an intermediate 
position between the zero-dimensional metrizable spaces and 
the strongly zero-dimensional metrizable spaces (a metrizable 
space $X$ has dimension $\dim$ zero if and only if it is metrizable 
by a non-Archimedean metric, and  this non-Archimedean metric can be 
assumed to take only rational values (see~\cite{Eng}). The Graev 
extension~\cite{Graev} of such a metric to the free group $F(X)$ 
takes only rational values as well; therefore, the group 
$F(X)$ with the Graev metric has dimension $\ind$ zero, and it contains 
$X$ as a subspace). In the second section, we construct 
a space (this is a special case of Mrowka's space $\mu\nu_0$) 
which can be embedded as a closed subspace in a 
zero-dimensional metrizable group 
but is not strongly zero-dimensional; 
thereby, an example 
of a metrizable group with noncoinciding dimensions $\ind$ and $\dim$ 
is obtained. The third section contains an example of a 
zero-dimensional metrizable 
space which cannot be embedded in a metrizable zero-dimensional 
group. 

\section{Spaces Embeddable in Zero-Dimensional Metrizable Groups}

The purpose of this section is to prove the following theorem. 

\begin{thm}\label{theorem1.1}
A topological space $X$ can be embedded in 
a metrizable topological group 
with dimension $\ind$ zero
if and only if the topology of $X$ is generated by 
a uniformity which 
has a countable base consisting of open-and-closed 
sets.
\end{thm}

The ``only if'' part is obvious: if $X$ is embedded in a group 
$G$ and clopen sets $U_n$, where $n \in \omega$, 
form a neighborhood base at the identity in $G$, then the required base 
of a uniformity on $X$ consists of the entourages 
$\mathcal U_n=\{(x, y): xy^{-1}\in U_n\cap U_n^{-1}\}$. 

The rest of this section is devoted to the proof of the 
reverse implication.  By $A(X)$ we denote the free Abelian group 
generated by $X$; the letters  $a$, $b$, $u$, $v$, $w$, $x$, $y$, and 
$z$ always denote elements of $X$, the letters $i$, $j$, $k$, $l$, 
$m$, $n$, $r$, $s$, $t$, and $N$ denote nonnegative integers, 
and $g$ and $h$ denote elements of the free Abelian group 
$A(X)$.  We use the definition of uniformities and entourages given 
in~\cite{Eng}; in particular, all entourages are assumed to be 
symmetric. For $A, B\subset X\times X$, we write 
$$%
A\circ B=\{(x,y)\in X\times X: 
\text{there exists a $z\in X$ for which $(x, z)\in A$ and $(z,y)\in 
B$}\}.
$$%
If $A$ or $B$ is a one-point set, we omit 
the braces in the notation of this set and write, e.g.,  
$A\circ (x,y)$. In particular, $(x, y)\circ (y, z) = (x, z)$ and 
$(x, y)\circ (u, z) =\varnothing$ if $y\ne u$. 

If $(x, y)= 
(x=x_1, y_1)\circ
(y_1=x_2, y_2)\circ
\dots \circ
(y_{n-1}=x_n, y_n=y)$,
then, obviously, 
$x-y=\sum_{i=1}^n(x_i-y_i)$ in $A(X)$. We write 
$$%
x-y=\dcircledsum_{i=1}^n(x_i-y_i)
$$%
in this case.

\begin{lem}\label{lemma0}
Let $\mathcal V_0, \mathcal V_1, 
\dots$ be \textup(symmetric\textup) elements 
of a uniformity of a set $X$ 
such that $\mathcal V_0=X\times X$ and $\mathcal V_{i+1}\circ
\mathcal V_{i+1}\circ\mathcal V_{i+1} \subset 
\mathcal V_i$ for $i=1, 2, \dots$\,, and let 
$\mathcal U_i = \mathcal V_{i^2}$ for $i\in \omega$. 
Suppose that $\{k_1, \dots, k_n\}$ is a set of positive integers 
in which each number $i$ occurs at most $i$ times.  Then 
$\mathcal U_{k_1}\circ \mathcal U_{k_2}\circ \dots \circ \mathcal 
U_{k_n}\subset \mathcal U_{k_*-1}$, where $k_*=\min_i\{k_i\}$.  
\end{lem}

\begin{proof}
If $k_*=1$, then the assertion holds trivially. Suppose that 
$k_*>1$, i.e., all $k_i$ are larger than~1.  Let $\rho$ be a 
pseudometric on $X$ such that $\mathcal V_i\subset 
\{(x,y):\rho(x,y)\le \frac1{2^i}\}\subset \mathcal V_{i-1}$ for any 
$i\ge 1$ (it exists by Theorem~8.1.10 from~\cite{Eng}). For 
$(x,y)\in \mathcal U_{k_1}\circ \mathcal U_{k_2}\circ \dots \circ 
\mathcal U_{k_n}$, we have 
$$%
(x,y)=(x=z_1, z_2)\circ(z_2, z_3)\circ \dots \circ (z_{n-1}, 
z_n)\circ (z_n, z_{n+1}=y), 
$$%
where 
$(z_i,  z_{i+1})\in \mathcal U_{k_i}=\mathcal V_{k_i^2}$ for $i\le n$. 
Hence  
$$%
\rho(x, y)\le \sum_{i=1}^n\frac1{2^{k^2_i}}\le 
\sum_{j=k_*}^\infty \frac j{2^{j^2}} \le 
\sum_{j=k_*}^\infty \frac {2^{j-1}}{2^{j^2}}
\le \sum_{j=k_*}^\infty \frac 1{2^{j^2-j+1}} =
\frac 1 {2^{k_*^2-k_*}}\le \frac 1 {2^{(k_*-1)^2+1}}.
$$%
Therefore, $(x,y)\in 
\mathcal V_{(\min_i\{k_i\}-1)^2}=
\mathcal U_{\min_i\{k_i\}-1}$.
\end{proof}

Let $X$ be a topological space whose topology is 
generated by a uniformity $\boldsymbol{\mathcal W}$ having 
a countable base $\{\mathcal W_n\}$ consisting of clopen 
sets.  Take a sequence $\mathcal V_0, \mathcal V_1, 
\dots$ of clopen entourages  
such that $\mathcal V_0=X\times X$, $\mathcal V_1=\mathcal W_1$, 
and $\mathcal V_{i+1}\circ\mathcal V_{i+1}\circ\mathcal V_{i+1} \subset 
\mathcal V_i\cap \mathcal W_{i+1}$ 
for $i= 2,3, \dots$\,. We set 
$\mathcal U_i = \mathcal V_{i^2}$ for $i\in \omega$. The sequence 
$\boldsymbol{\mathcal U}=\{\mathcal U_i\}$ 
is a 
base of the uniformity $\boldsymbol{\mathcal W}$,  
and the sets 
$$%
W_n(\boldsymbol{\mathcal U}) = 
\bigcup_{k\in \omega} \bigl\{\sum_{i=1}^k(x_i-y_i):
(x_i, y_i)\in \mathcal U_{n\cdot i}\bigr\}
$$%
form a neighborhood base at zero for  some group topology 
$\mathcal T_{\boldsymbol{\mathcal U}}$ on the free Abelian
group $A(X)$ which  
induces the initial topology (generated by the uniformity 
$\boldsymbol{\mathcal W}$) on $X$.  Indeed, 
it is easy to show that $2W_{2n}(\boldsymbol{\mathcal U})
\subset W_n(\boldsymbol{\mathcal U})$ for $n\ge 1$ 
and that if $g=\sum_{i=1}^k(x_i-y_i)\in 
W_n(\boldsymbol{\mathcal U})$, 
then $g+ W_{n(k+1)}(\boldsymbol{\mathcal U})
\subset W_n(\boldsymbol{\mathcal U})$; in addition, all 
sets $W_n(\boldsymbol{\mathcal U})$ are symmetric and 
contain the empty word (the zero of the group $A(X)$),  
and $W_n(\boldsymbol{\mathcal U})\cap W_k(\boldsymbol{\mathcal U})\supset 
W_{\max\{k,n\}}(\boldsymbol{\mathcal U})$. To see that 
$\mathcal T_{\boldsymbol{\mathcal U}}$ 
induces the topology generated by the uniformity 
$\boldsymbol{\mathcal W}$ with base $\{\mathcal W_i\}$ on $X$, 
it suffices to note that, for any $x\in X$ and $n\ge 1$, we have 
\begin{multline*}
(x+W_n(\boldsymbol{\mathcal U}))\cap X=
\{x+(y-x): y-x\in W_n(\boldsymbol{\mathcal U})\}\\
=\{y\in X: (y, x)\in \bigcup\{U_{n\cdot \pi(1)}
\circ \dots\circ U_{n\cdot \pi(k)}: k\ge 1, \ \pi\in S_k\}\}
\end{multline*}
(here $S_k$ is the permutation group on $\{1, \dots, k\}$). 
By Lemma~\ref{lemma0}, 
$$%
(x+W_n(\boldsymbol{\mathcal U}))\cap X\subset 
\{y\in X: (y, x)\in U_{n-1}\}.
$$%
On the other hand, clearly, 
$$%
(x+W_n(\boldsymbol{\mathcal U}))\cap X\supset 
\{y\in X: (y, x)\in U_{n}\}.
$$%

Our immediate goal is to construct 
a base of the topology $\mathcal T_{\boldsymbol{\mathcal U}}$ on 
$A(X)$ consisting of open-and-closed (in this topology) sets.  

\begin{defn}\label{defd}
For $x, y\in X$, we set 
$$%
d(x, y) = \begin{cases} 
\frac 1{\max\{k: (x, y)\in \mathcal U_k\}}
&\text{if $x\ne y$},\\
0&\text{if $x= y$}.
\end{cases}
$$%
Thus, for $x\ne y$, the number 
$d(x,y)$ is uniquely determined by the conditions 
$(x, y)\in \mathcal U_{\frac 1{d(x, y)}}$ and 
$(x, y)\notin \mathcal U_{\frac 1{d(x, y)}+1}$.
\end{defn}

\begin{defn}\label{def*}
Suppose that $k\in \omega$, $x_i, y_i\in X$ for $i\le k$, and 
$$%
g=\sum_{i=1}^k (x_i-y_i)\in A(X).
$$%
We say that the sum (decomposition) 
$\sum_{i=1}^k (x_i-y_i)$  \emph{satisfies condition $(*)$} 
if 
$$%
d(x_i, y_j)\ge \min\{d(x_i, y_i), d(x_j, y_j)\} 
\quad\text{for any $i, j\le k$.}\eqno{(*)}
$$%
Sometimes, when it is clear what decomposition of $g$ is meant, we 
say the word $g$ itself satisfies condition $(*)$ 
(meaning that condition $(*)$ holds for the decomposition).  
\end{defn}

\begin{rem}\label{remark1}
Suppose that $d(x,y)\le d(x, y_i)$ and $d(x,y)\le d(x_i, y)$ 
for all $i\le k$. 
Then $\sum_{i=1}^k (x_i-y_i)$ satisfies condition $(*)$ 
if and only if 
$\sum_{i=1}^k (x_i-y_i)+(x-y)$ satisfies condition $(*)$. Moreover, 
if $\sum_{i=1}^k (x_i-y_i)$ satisfies condition $(*)$, 
then  $\sum_{i\in I} (x_i-y_i)$ satisfies condition $(*)$ for any 
$I\subset \{1, \dots, k\}$.
\end{rem}

\begin{lem}\label{lemma1}
Suppose that 
\begin{enumerate}
\item[\upshape (1)]
$g=\sum_{i=1}^k(x_i-y_i)$\textup; 
\item[\upshape (2)]
$(x_i, y_i)= 
(x_i=x_i^{(1)}, y_i^{(1)})\circ
(y_i^{(1)}=x_i^{(2)}, y_i^{(2)})\circ
\dots \circ
(y_i^{(k_i-1)}=x_i^{(k_i)}, y_i^{(k_i)}=y_i)$, i.e., 
$x_i-y_i=\circledsum_{j=1}^{k_i}(x_i^{(j)}-y_i^{(j)})$ for 
each $i\le k$\textup; 
\item[\upshape (3)] 
$(x_i^{(j)},y_i^{(j)})\in 
\mathcal U_{N\cdot n_i^{(j)}}$ for all $i\le k$ and $j\le 
k_i$\textup; 
\item[\upshape (4)] 
if $m\le k$, then $n_i^{(j)}=m$ for 
at most one pair $i, j$\textup; 
\item[\upshape (5)] 
if $m> k$, then 
$n_i^{(j)}=m$ for at most $m-k+1$ pairs $i, j$\textup; 
\item[\upshape (6)] $k> 1$.  
\end{enumerate} 
Then $g=\sum_{i=1}^{k-1}({x'}_i-{y'}_i) 
+ x'' -y''$, where each of the letters ${x'}_i$, ${y'}_i$, $x''$ and 
$y''$ is contained in one of the decompositions from \textup{(2)} and 
$\sum_{i=1}^{k-1}({x'}_i-{y'}_i)$ satisfies conditions 
\textup{(2)--(5)} with $x_i$ replaced by ${x'}_i$, $x_i^{(j)}$ by
${x'}_i^{(j)}$, $y_i$ by ${y'}_i$, $y_i^{(j)}$ by ${y'}_i^{(j)}$, 
$k$ by $k-1$, $k_i$ by $k'_i$, and $n_i^{(j)}$ by ${n'}_i^{(j)}$\textup; 
moreover, 
\begin{enumerate} 
\item[\upshape (7)] $d(x'', y'')\le d(x'', 
{y'}_i^{(j)})$ and $d(x'', y'')\le d({x'}_i^{(j)}, y'')$  for all 
$i\le k-1$ and $j\le k'_i$\textup; 
\item[\upshape (8)] $(x'', y'') \in \mathcal U_{N\cdot k-1}$.  
\end{enumerate} 
\end{lem}

\begin{proof}
Take any pair $(u, v)$ for which 
$u\in \{x_1, \dots, x_k\}$, $v\in \{y_1, \dots, y_k\}$ and 
$d(u, v)$ is minimal; if there exists a pair of the form 
$(x_i, y_i)$ with these 
properties, then let $(u, v)$ be such a pair. By condition (4), 
there exists an $i\le k$ for which $\min_{j\le k_i}n_i^{(j)}\ge k$.  
Conditions (2)--(5) and Lemma~\ref{lemma0} imply 
that $(x_i, y_i)\in \mathcal U_{N\cdot k-1}$ for this $i$, i.e., 
$d(x_i, y_i)\le \frac 1{N\cdot k-1}$.  Therefore, $d(u,v)\le \frac 
1{N\cdot k-1}$ (by virtue of minimality), i.e., $(u, v)\in \mathcal 
U_{N\cdot k-1}$.  If $(u,v)=(x_i, y_i)$ for some $i\le k$, then 
the required decomposition consists 
of the term $x_i - y_i$ and the sum of all other terms in 
the initial decomposition of the word $g$; in other words, 
it suffices to set ${x'}_j = x_j$  and ${y'}_j = y_j$ for $j<i$, 
${x'}_j = x_{j+1}$ and ${y'}_j = y_{j+1}$ for $j=i, \dots, k-1$, 
$x''=x_i$, and $y''= y_i$.  The decompositions from (2) remain the 
same for all $x'_j-y'_j$. 

If $u=x_i$, $v= y_j$, and $i\ne j$, i.e., the function $d$ does not 
attains its minimum for pairs of the form $(x_r, y_r)$, 
then $d(u, v)<\frac 
1{N\cdot k-1}$, because, as mentioned above, $d(x_s, y_s)\le \frac 
1{N\cdot k-1}$ for some $s$. Therefore, $d(u, v) \le \frac1{N\cdot 
k}$. Without loss of generality, we can assume that $i< j$. 
We set ${x'}_r=x_r$ and ${y'}_r=y_r$ for $r< j$ such that $r\ne i$,  
${x'}_i = x_j$, $y'_i = y_i$, $x'_r= x_{r+1}$ and $y'_r= y_{r+1}$ for 
$r=j, \dots, k-1$,  
$x''=x_i$, and $y'' = y_j$; in 
other words, we replace the pairs $(x_i, y_i)$ and $(x_j, y_j)$ by 
$(x_j, y_i)$ and $(x_i, y_j)$.  The fulfillment of condition (7) 
follows from the choice of the pair $(u, v)$, and (8) holds because 
$d(u, v)<\frac 1{N\cdot k-1}$. 

The decompositions from (2) and numbers of the form $n_r^{(t)}$ 
remain the same for the pairs $({x'}_r, {y'}_r)=(x_r, y_r)$ with 
$r\ne i$, which coincide with $(x_r, y_r)$ or  $(x_{r+1}, y_{r+1})$; 
for $({x'}_i, {y'}_i)=(x_j, y_i)$, we set 
$k'_i= k_i+k_j+1$ and take the decomposition 
$$%
{x'}_i-{y'}_i=x_j-y_i= x_j-y_j+y_j-x_i+x_i-y_i
=
\dcircledsum_{r=1}^{k_j}(x_j^{(r)}-y_j^{(r)}) + (y_j-x_i)+
\dcircledsum_{s=1}^{k_i}(x_i^{(s)}-y_i^{(s)});
$$%
thus, we set 
${x'}_i^{(t)}=x_j^{(t)}$ and ${y'}_i^{(t)}=y_j^{(t)}$ for $t\le k_j$, 
${x'}_i^{(k_j+1)}=y_j$, ${y'}_i^{(k_j+1)}=x_i$, 
${x'}_i^{(t)}=x_i^{(t-k_j-1)}$, and ${y'}_i^{(t)}=y_i^{(t-k_j-1)}$ for 
$t=k_j+2, \dots, k_j+k_i+1$.
As mentioned above, $(u, v)=(y_j, x_i)\in \mathcal U_{N\cdot k}$. 
Therefore, setting ${n'}_i^{(t)}=n_j^{(t)}$ for $t\le k_j$ and 
${n'}_i^{(k_j+1)}=k$ 
and ${n'}_i^{(t)}=n_i^{(t-k_j-1)}$ for 
$t=k_j+2, \dots, k_j+k_i+1$, we obtain 
$({x'}_i^{(t)},{y'}_i^{(t)})\in \mathcal U_{N\cdot {n'}_i^{(t)}}$ 
for all $t\le 
k'_i$. The term ${x'}_i^{(k_j+1)}-{y'}_i^{(k_j+1)}=y_j-x_i$ is the 
only new element in the sum 
$$%
\sum_{r=1}^{k-1}\dcircledsum_{s=1}^{k'_r}({x'}_r^{(s)}-{y'}_r^{(s)})
=\sum_{r=1}^{k-1}({x'}_r-{y'}_r)
$$%
in comparison with the sum 
$$%
\sum_{r=1}^{k}\dcircledsum_{s=1}^{k_r}(x_r^{(s)}-y_r^{(s)})
=\sum_{r=1}^{k}(x_r-y_r),
$$%
and we have ${n'}_i^{(k_j+1)}=k>k-1$ for this element; 
the numbers of the form ${n'}_r^{(s)}$ corresponding to the other 
terms are equal to the numbers corresponding to them as terms of 
the sum $\sum_{r=1}^{k}\circledsum_{s=1}^{k_r}(x_r^{(s)}-y_r^{(s)})$.  
Therefore, 
$\sum_{r=1}^{k-1}({x'}_r-{y'}_r)
=\sum_{r=1}^{k-1}\circledsum_{s=1}^{k'_r}({x'}_r^{(s)}-{y'}_r^{(s)})$ 
satisfies condition (4) 
with $k$ replaced by $k-1$ and $n_i^{(j)}$ by ${n'}_i^{(j)}$; 
it also satisfies the part of condition (5) (with 
the appropriate replacements) that relates to the 
number of ${n'}_r^{(s)}>k$.  By condition (4), the sum 
$\sum_{r=1}^{k}\circledsum_{s=1}^{k_r}(x_r^{(s)}-y_r^{(s)})$ contains at 
most one term for which $n_r^{(s)}=k$.  Therefore, the sum 
$\sum_{r=1}^{k-1}\circledsum_{s=1}^{k'_r}({x'}_r^{(s)}-{y'}_r^{(s)})$ 
contains at most two terms for which ${n'}_r^{(s)}=k$; thus, 
condition (5) 
with $k$ replaced by $k-1$ and $n_i^{(j)}$ by ${n'}_i^{(j)}$ is 
satisfied fully.  Conditions (2) and (3) with 
the appropriate replacements hold by construction.  
\end{proof}

\begin{cor}\label{corollary1}
If $k,n\in \omega$, 
$g=\sum_{i=1}^k(x_i-y_i)$, 
and $(x_i,y_i)\in \mathcal U_{(n+1)\cdot i}$ for all $i\le k$, 
then 
$g=\sum_{i=1}^k(\tilde x_i-\tilde y_i)$, 
where $(\tilde x_i,\tilde y_i)\in \mathcal U_{n\cdot i}$ 
for all $i\le k$ 
and the decomposition $\sum_{i=1}^k(\tilde x_i-\tilde y_i)$ 
satisfies condition 
$(*)$.
\end{cor}

\begin{proof}
This assertion is proved by repeatedly applying 
Lemma~\ref{lemma1}  with $N= n+1$ to the word $g$. 
If $k\le1$, then the assertion holds trivially.  
If $k>1$, then we can apply Lemma~\ref{lemma1} with $N= n+1$ 
and obtain a decomposition 
$$%
g=\sum_{i=1}^{k-1}({x'}_i-{y'}_i) + x'' -y'' 
$$%
with the properties described in the lemma. 
We have $(x'',y'')\in \mathcal U_{N\cdot k -1} 
\subset \mathcal U_{n\cdot k}$ and, for each $i\le k-1$,
${x'}_i-{y'}_i=\circledsum_{j=1}^{k'_i}({x'}_i^{(j)}-{y'}_i^{(j)})$, where 
$({x'}_i^{(j)},{y'}_i^{(j)})\in \mathcal U_{(n+1)\cdot 
{n'}_i^{(j)}}$; moreover, if $m>k-1$, then ${n'}_i^{(j)}=m$ for 
at most $m-k+2$ pairs $i, j$, and if $m\le k-1$, then 
${n'}_i^{(j)}=m$ for at most one pair $i, j$. We apply 
Lemma~\ref{lemma1} first to the sum 
$\sum_{i=1}^{k-1}({x'}_i-{y'}_i)$, then to the obtained 
decomposition, then to the new decomposition,  
and so on, while possible; in the end (after $k-1$ steps), 
we obtain a decomposition
$$%
g=\dcircledsum_{j=1}^{\tilde k'_1}(\tilde {x'}_1^{(j)}-\tilde {y'}_1^{(j)}) + 
\sum (\tilde x'' - \tilde y''),
$$%
where 
$\circledsum_{j=1}^{\tilde k'_i}(\tilde {x'}_1^{(j)}-\tilde {y'}_1^{(j)})=
\tilde x-\tilde y$ 
for some $\tilde x, \tilde y\in X$ and 
$\sum (\tilde x'' - \tilde y'')$ 
denotes the sum of the residual terms of the form $x''-y''$ 
obtained at all steps. The pairs of letters in each residual 
term belongs to the entourage $\mathcal U_{N\cdot 
(k-s+1)-1}\subset U_{n\cdot (k-s+1)}$, where $s<k$ is the  
number of the step at which this term has appeared (and 
$N=n+1$).  Moreover, 
$$%
(\tilde x, \tilde y)=
(\tilde {x'}_1^{(1)},\tilde {y'}_1^{(1)})\circ\dots \circ
(\tilde {x'}_1^{(\tilde k'_1)},{y'}_1^{(\tilde k'_1)})
 \in 
\mathcal U_{N\cdot \tilde {n'}_1^{(1)}}\circ \dots \circ 
U_{N\cdot \tilde {n'}_1^{(\tilde k'_1)}},
$$%
and, for $m=1, 2, \dots\,$, 
${n'}_1{(j)}=m$ for at most $m$ indices $j$. 
Therefore, by Lemma~\ref{lemma0}, 
$$%
(\tilde x, \tilde y)\in 
\mathcal U_{N\cdot \min_{j\le \tilde k'_1}\{\tilde {n'}_1^{(j)}\}-1}\subset 
\mathcal U_{N-1}= \mathcal U_{n}.
$$%
Condition (7) from lemma~\ref{lemma1} and 
Remark~\ref{remark1}, as well as the fact that no new letters appear 
in repeatedly applying Lemma~\ref{lemma1},
ensure the fulfillment of condition $(*)$.  
\end{proof}

\begin{lem}\label{lemma2}
Suppose that $I=\{k_1, \dots, k_l\}$ is a finite set of 
different positive integers enumerated in increasing order, 
$g=\sum_{i=1}^k (x_i-y_i)$, $h=\sum_{j=1}^l (u_j-v_j)$, the decompositions 
$\sum_{i=1}^k (x_i-y_i)$ and $\sum_{j=1}^l (u_j-v_j)$ satisfy 
condition $(*)$, $(x_i, y_i)\in \mathcal U_{n\cdot i}$ for $i\le k$, 
and $(u_j, v_j)\in \mathcal U_{k_j}$ for $j\le l$\textup; suppose 
also that if $i,j\le k$ and $F=(f_1, \dots, f_r)$, $F'=(f'_1, \dots, 
f'_{r'})$ are finite ordered sequences of elements of $I\cup 
\{k_l+1, k_l+2, \dots\}$ in each of which every 
element of $I$ occurs at most once and every 
positive integer $s>k_l$ occurs at most $s$ times, then 
\begin{enumerate} 
 \item[(1)] 
 $\mathcal U_{f_1}\circ \dots 
\circ\mathcal U_{f_r}\circ (x_i, y_j)\circ \mathcal U_{f'_1}\circ 
\dots \circ\mathcal U_{f_{r'}} \subset \mathcal U_{\frac1{d(x_i, 
y_j)}}$ 
\quad 
and 
\item[(2)] $\mathcal U_{f_1}\circ \dots 
\circ\mathcal U_{f_r}\circ (x_i, y_j)\circ \mathcal U_{f'_1}\circ 
\dots \circ\mathcal U_{f'_{r'}} \cap \mathcal U_{\frac1{d(x_i, 
y_j)}+1}=\varnothing$.  
\end{enumerate} 
Then $g+h=\sum_{i=1}^m 
(z_i-w_i)$, where $m\le k+l$, the decomposition 
$\sum_{i=1}^m (z_i-w_i)$ satisfies 
condition $(*)$, $z_i\in \{x_1, \dots, x_k, u_1, \dots, u_l\}$, 
$w_i\in \{y_1, \dots, y_k, v_1, \dots, v_l\}$, $(z_i, w_i)\in 
 \mathcal U_{n\cdot i}$ for $i\le k$, and $(z_{k+i}, w_{k+i})\in 
 \mathcal U_{k\cdot i}$ for $i\le m-k$ (if $m>k$).  
\end{lem}

\begin{proof}
First, note that (2) implies $d(x_i, y_j) > \frac 1{k_1}\ge \frac 1{k_s}$ 
for any $i,j\le k$ and $s\le l$. Indeed, otherwise,  
$\mathcal U_{k_1}\circ (x_i, y_j)\ni (y_j, y_j)$; clearly, 
$(y_j, y_j)\in \mathcal U_{\frac 1{d(x_i, y_j)}+1}$, while by condition~(2), 
$\mathcal U_{k_1}\circ (x_i, y_j)\cap \mathcal U_{\frac1{d(x_i, 
y_j)}+1}=\varnothing$ (consider $F=\{k_1\}$ and $F'=\varnothing$). 
This implies, in particular, that $k_1> n\cdot k$.

We shall prove the lemma by induction on $l$. 
If $l=0$ (i.e., the word $h$ is empty), 
then the assertion holds trivially.  Suppose that $l>0$ and 
the assertion is true for smaller $l$. Choose 
$h'\in \{x_1, \dots, x_k, u_1, \dots, u_l\}$ and $h''\in \{y_1, 
\dots, y_k, v_1, \dots, v_l\}$ for which $d(h', h'')$ is minimal. 
Since $h$ is nonempty, we have 
\begin{itemize} 
\item[(i)] 
$d(h', h'')\le \frac 1{k_l}$ (because $d(u_l, v_l)\le \frac 1{k_l}$ and 
$d(h', h'')$ is minimal) and 
\item[(ii)] 
either $h'\in \{u_1, \dots, 
u_l\}$ or $h''\in \{v_1, \dots, v_l\}$ (this follows from 
(i) and because, by condition (2), $d(x_i, y_j)> \frac 1{k_l}$ for 
all $i, j\le k$); moreover, we can assume 
that if $h'= u_i$ and 
$h''= v_j$, then $i=j$; otherwise, we replace the pair $h', h''$ 
by the pair $u_i, v_i$ or $u_j, v_j$ for which the value 
of $d$ does not exceed $d(h', h'')$ (such a pair exists 
because the decomposition $\sum_{i=1}^l (u_i-v_i)$ satisfies condition 
$(*)$).  
\end{itemize} 
If $h'= u_i$ and $h''=v_i$ for some $i\le l$, 
then we set $\tilde u_j = u_j$ and $\tilde v_j = v_j$ for $j< i$,  
$\tilde u_j = u_{j+1}$ and $\tilde v_j = v_{j+1}$ 
for $j=i, \dots, l-1$, 
$\tilde h=\sum_{i=1}^{l-1}(\tilde u_i-\tilde v_i)=h-(u_i-v_i)$, 
and 
$\tilde I=\{k_1, \dots, k_{l-1}\}$.
Note that the conditions of the 
lemma hold for $\tilde I$, 
$g=\sum_{i=1}^k (x_i-y_i)$, and $\tilde h$. 
By the induction assumption, $g+\tilde h=\sum_{i=1}^{\tilde m} 
(\tilde z_i-\tilde w_i)$, 
where 
$\tilde m\le k+l-1$, the decomposition 
$\sum_{i=1}^m (\tilde z_i-\tilde w_i)$ satisfies 
condition $(*)$, 
$\tilde z_i\in \{x_1, \dots, x_k, \tilde u_1, \dots, \tilde u_{l-1}\} 
=\{x_1, \dots, x_k, u_1, \dots, u_l\}\setminus \{u_i\}$, 
$\tilde w_i\in \{y_1, \dots, y_k, \tilde v_1, \dots, \tilde v_{l-1}\}
=\{x_1, \dots, x_k, v_1, \dots, v_l\}\setminus \{v_i\}$, 
$(\tilde z_i, \tilde w_i)\in 
 \mathcal U_{n\cdot i}$ for $i\le k$, and $(\tilde z_{k+i}, 
\tilde w_{k+i})\in 
 \mathcal U_{k\cdot i}$ for $i\le \tilde m-k$ (if $\tilde m>k$).
By the definition of the pair $h'=u_i$, $h''=v_i$ and Remark~\ref{remark1}, 
the decomposition $g+ h=\sum_{i=1}^{\tilde m} 
(\tilde z_i-\tilde w_i)+(h'-h'')$
has the required properties (recall that 
$d(h', h'')\le \frac 1{k_l}<\frac 1{n\cdot k}$).  

Suppose that 
$h'$ and $h''$ cannot be chosen among the letters 
of the form $u_i$ and $v_j$, i.e., either $h'= x_i$ and $h'' = v_j$ 
for some $i \le k$ and $j\le l$ and $d(x_i, v_j)< d(u_r, v_s)$ for 
all  $r, s\le l$ (i.e., $d(x_i, v_j)< \frac1{k_l}$) or $h'= 
u_i$ and $h'' = y_j$ for some $i \le l$ and $j\le k$ and $d(u_i, 
y_j)< d(u_r, v_s)$ for all  $r, s\le l$ (i.e., $d(u_i, y_j)< 
\frac1{k_l}$).  For definiteness, suppose that $h'= x_i$ and $h'' = 
v_j$.  We have $(u_j, x_i)= (u_j, v_j)\circ (v_j, x_i)\in \mathcal 
U_{k_j}\circ \mathcal U_{k_l}$, and conditions (1) and (2) 
imply $d(u_j, y_r)=d(x_i, y_r)$ for all $r\le k$.  We set 
$\tilde x_i=u_j$, $\tilde x_s=x_s$ for $s\ne i$, and $\tilde y_r=y_r$ 
for all $r\le k$; thus, the word $\sum_{s=1}^k (\tilde x_s -\tilde y_s)$ 
differs from $\sum_{s=1}^k (x_s - y_s)$ in one letter $\tilde 
x_i$, and $d(\tilde x_s, \tilde y_t)=d(x_s, y_t)$ for all $s,t\le k$ 
(this means that $\sum_{s=1}^k(\tilde x_s, \tilde y_t)$ satisfies 
condition $(*)$ and $(\tilde x_i, \tilde y_i)\in \mathcal U_{n\cdot 
i}$ for $i\le k$).  We also set $\tilde u_s=u_s$ and $\tilde 
v_s=v_s$ for $s<j$, $\tilde u_s=u_{s+1}$ and $\tilde v_s=v_{s+1}$ for 
$s=j, \dots, l-1$, $\tilde u_l=x_i$, and $\tilde v_l = v_j$; thus, 
the word $\sum_{s=1}^l (\tilde u_s - \tilde y_s)$ is obtained from 
$\sum_{s=1}^l (u_s - y_s)$ by deleting the term $u_j-v_j$ and 
inserting $\tilde u_l-\tilde v_l= x_l-v_l$.  Since 
$\tilde u_l=x_i$, $\tilde v_l=v_j$, and $d(x_i, v_j)$ is 
minimal, it follows that $d(\tilde u_l, \tilde v_l)\le 
d(\tilde u_l, \tilde v_r)$ and $d(\tilde u_l, \tilde v_l)\le d(\tilde 
u_r, \tilde v_l)$ for all $r< l$.  Therefore, the word $\sum_{s=1}^l 
(\tilde u_s - \tilde v_s)$ satisfies condition $(*)$.  Indeed, the word 
$h$ satisfies condition $(*)$; according to Remark~\ref{remark1}, 
deleting the term $u_j-v_j$ does not violate condition $(*)$; 
applying Remark~\ref{remark1} again with taking into account 
the minimality of $d(\tilde u_l, \tilde v_l)$, we conclude that 
$\sum_{s=1}^l (\tilde u_s - \tilde v_s)$ satisfies condition $(*)$.  
We set $\tilde k_s=k_s$ for $s= 1, \dots, j-1$, $\tilde 
k_s=k_{s+1}$ for $s=j,\dots, l-1$, $\tilde k_l=k_l+1$, and 
$$%
\tilde I=\{\tilde k_1, \dots, \tilde k_l\} = 
(I\setminus \{k_j\})\cup \{k_{l}+1\}.  
$$%
We have 
$(\tilde u_s, \tilde v_s)\in \mathcal U_{\tilde k_s}$ for $s\le l$. 
Finally, $\tilde k_l=k_l+1$ and $\tilde I$ does not contain $k_j$; 
therefore, if $F=(f_1, \dots, f_r)$ is a finite ordered sequence 
of elements of the set $\tilde I\cup \{\tilde k_l+1, \tilde k_l+2, \dots\}$ 
with the properties (a)~each element from $\tilde I$ occurs 
in $F$ at most once and (b)~each element $s$ larger than all 
elements of $\tilde I$ 
occurs at most $s$ times, then 
the sequences $F$ and $(f_1, \dots, f_r, k_j, k_l+1)$ have the  
same properties with respect to the set $I$. This observation, 
conditions (1) and (2) of the lemma being proved, and the 
relations 
$$%
(\tilde x_i, \tilde y_t)=
(u_j, y_t)=(u_j, v_j)\circ(v_j, x_i)\circ(x_i, y_t)\in 
\mathcal U_{k_j}\circ\mathcal U_{k_l+1}\circ(x_i, y_t)
$$%
and $(\tilde x_s, \tilde y_t)=(x_s, y_t)$ for $s\ne i$ and any $t$ 
imply that, for any $s,t\le k$ and any two 
finite ordered sequences
$(f_1, \dots, f_r)$ and $(f'_1, \dots, f'_{r'})$ of elements 
of the set $\tilde I\cup \{\tilde k_l+1, \tilde k_l+2, \dots\}$ in 
each of which
every element of $\tilde I$ occurs at most once and 
every element $s>\tilde k_l$ occurs at most $s$ times, 
we have 
\begin{enumerate} 
\item[($\tilde 1$)] $\mathcal 
U_{f_1}\circ \dots \circ\mathcal U_{f_r}
\circ (\tilde x_s, 
\tilde y_t)\circ \mathcal U_{f'_1}\circ \dots \circ\mathcal 
U_{f'_{r'}} \subset \mathcal U_{\frac1{d(\tilde x_s, 
\tilde y_t)}}$
and
\item[($\tilde 2$)]
$\mathcal U_{f_1}\circ \dots\circ \mathcal 
U_{f_r}\circ (\tilde x_s, \tilde y_t)\circ \mathcal 
U_{f'_1}\circ \dots \circ\mathcal U_{f'_{r'}}
\cap \mathcal U_{\frac1{d(\tilde x_s, \tilde y_t)}+1}=\varnothing$.
\end{enumerate} 
Thus, the set $\tilde I$ and the words 
$\sum_{s=1}^k (\tilde x_s - \tilde y_s)$
and 
$\sum_{s=1}^l (\tilde u_s - \tilde v_s)$
satisfy the conditions of the  lemma. Moreover, the set of letters 
(with signs) of which these words consist 
coincides with the set of letters in the words 
$\sum_{s=1}^k (x_s - y_s)$ and $\sum_{s=1}^l (u_s - v_s)$; therefore, 
the function $d$ takes minimal value at the same pair of letters 
$(h', h'') = (x_i, v_j)=(\tilde u_j, \tilde v_j)$.  However, these 
letters form a summand in the decomposition $\sum_{s=1}^l 
(\tilde u_s - \tilde v_s)$; this situation was considered 
at the beginning of the proof. As there, we delete this summand,  
apply the induction assumption, and insert the deleted summand back; 
as a result, we obtain a representation 
$g+ h=\sum_{i=1}^{\tilde m} 
(\tilde z_i-\tilde w_i)+(h'-h'')$, where 
$\tilde m\le k+l-1$, 
$(\tilde z_i, \tilde w_i)\in 
 \mathcal U_{n\cdot i}$ for $i\le k$, $(\tilde z_{k+i}, 
\tilde w_{k+i})\in 
 \mathcal U_{\tilde k\cdot i}$ for $i\le \tilde m-k$ (if $\tilde m>k$), 
and $(h',h'')\in \mathcal U_{\tilde k_l}$. Since $\tilde k_l=k_l+1\ge k_l$ 
and $\tilde k_i\ge k_i$ for all $i\le l-1$, this representation is as 
required.
\end{proof}

\begin{cor}\label{corollary2}
Suppose that $g=\sum_{i=1}^k (x_i-y_i)$, $h=\sum_{i=1}^l (u_i-v_i)$, 
the decompositions 
$\sum_{i=1}^k (x_i-y_i)$ and $\sum_{i=1}^l (u_i-v_i)$ 
satisfy condition $(*)$, 
$(x_i, y_i)\in \mathcal U_{n\cdot i}$ for $i\le k$, 
$(u_i, v_i)\in \mathcal U_{(N+1)\cdot i}$ for $i\le l$, 
$N\ge 2 n\cdot k$,
and, 
for any $i, j\le k$, 
\begin{enumerate} 
\item[(1)] $\mathcal U_N\circ (x_i, 
y_j)\circ \mathcal U_N \subset \mathcal U_{\frac1{d(x_i, y_j)}}$ and 
\item[(2)]
$\mathcal U_N\circ (x_i, y_j)\circ \mathcal 
U_N 
\cap \mathcal U_{\frac1{d(x_i, y_j)}+1}=\varnothing$.
\end{enumerate} 
Then $g+h=\sum_{i=1}^m (z_i-w_i)$, where the decomposition $\sum_{i=1}^m 
(z_i-w_i)$ satisfies condition $(*)$, 
$z_i\in \{x_1, \dots, x_k, u_1, \dots, u_l\}$, 
$w_i\in \{y_1, \dots, y_k, v_1, \dots, v_l\}$, 
and $(z_i, w_i)\in \mathcal 
U_{n\cdot i}$ for $i\le m$.
\end{cor}

\begin{proof}
This assertion follows immediately from Lemmas~\ref{lemma0} 
and~\ref{lemma2}.  
\end{proof}

\begin{lem}\label{lemma3}
Suppose that $n\ge 2$ and 
\begin{enumerate}
\item[(1)]
$\sum_{i=1}^n(z_i-w_i)$ satisfies condition $(*)$\textup; 
\item[(2)]
$\sum_{i=1}^{n-1}(z_{i+1}-w_i)$ satisfies condition $(*)$\textup; 
\item[(3)] 
$(z_{i+1}, w_i)\in \mathcal U_{k_i}$ for all $i\le n-1$, and the 
$k_i$ are different 
positive integers larger than~\textup{1;} 
\item[(4)] 
$d(z_1, w_1)\ge d(z_i, w_i)$ and $d(z_n, w_n)\ge d(z_i, w_i)$ 
for $i= 2, \dots, n-2$\textup; 
\item[(5)] 
$k_m=\min\{k_1, \dots, k_{n-1}\}$\textup; 
\item[(6)] either\quad 
\textup{($6_{\text{left}}$)}~$d(z_1, w_1)\ge d(z_n, w_n)$\quad 
or\quad
\textup{($6_{\text{right}}$)}~$d(z_n, w_n)\ge d(z_1, w_1)$
\end{enumerate}
\textup(the last condition is included for convenience\textup).
Then there exists a one-to-one map 
$$%
f\colon \{2, \dots, n-1\}\to \{k_1, \dots, k_{n-1}\}\setminus \{k_m\}
$$%
such that 
$d(z_i, w_i)\le \frac1{f(i)-1}$ for $i=2, \dots, n-1$ and 
$d(z_n, w_n)\le \frac1{k_m-1}$ \textup(if 
$(6_{\mathrm{left}})$  holds\textup) or $d(z_1, w_1)\le \frac1{k_m-1}$ 
\textup(if $(6_{\mathrm{right}})$ holds\textup).
\end{lem}

\begin{proof}
We prove the lemma by induction. For $n=2$, 
the map $f$ is trivial, 
and $d(w_1, z_2)\ge \min\{d(w_1, z_1), d(w_2, 
z_2)\}$ by condition $(*)$.  
This implies the required assertion, because it follows from 
(3) and (5) that $d(w_1, z_2) \ge \frac 1{k_m}$.  Suppose 
that $n>2$ and the assertion is true for smaller $n$. Let 
$m_1, \dots, m_r$ be the indices (or index) from $\{2, 
\dots, n-1\}$ for which the numbers $d(z_{m_j}, w_{m_j})$ 
are maximal (and equal to each other).  These indices divide the 
set of all indices into intervals.  Suppose that $m$ belongs to 
the $s$th interval, i.e., $k_i$ is minimal for $i\in \{m_s, 
m_s+1, \dots, m_{s+1}-1\}$, where $s=0, \dots, r$ (we assume that 
$m_0=1$ and $m_{r+1}=n$).  Suppose that $s>0$; for $s=0$, the 
argument is the same except that we must replace the 
conditions $j<s$ and $j\ge s$ by $j\le s$ and $j>s$ (that is, by 
$j=0$ and $j>0$), respectively, every time they are 
encountered.  Consider the words 
$z_{m_j}-w_{m_j}+z_{m_j+1}-w_{m_j+1}+ \dots 
+z_{m_{j+1}-1}-w_{m_{j+1}-1}+z_{m_{j+1}}-w_{m_{j+1}}$.  They
satisfy condition $(*)$, being subsums of a sum 
satisfying condition $(*)$, and to these words the induction 
hypothesis applies.  
Using the left version of the lemma for $j<s$ and the right version 
for $j\ge s$ (recall that, for $s=0$, the condition $j<s$ should be 
replaced by $j=0$ and the condition $j\ge s$, by $j>0$; in 
the situation under consideration, this means that if $s=0$, 
then the left version should be applied to $j=0$ and the right 
version, to $j>0$), 
we obtain one-to-one maps 
$$%
f_j\colon \{m_j+1, \dots, m_{j+1}-1\}\to 
\{k_{m_j}, \dots, k_{m_{j+1}-1}\}\setminus \{k_{\min\nolimits_j}\}
$$%
such that 
$d(z_{m_j+i}, w_{m_j+i})\le \frac1{f_j(m_j+i)-1}$ for $i=1, \dots, m_{j+1}
-m_j-1$ and 
$d(z_{m_{j+1}}, w_{m_{j+1}})\le \frac1{k_{\min\nolimits_j}-1}$ 
(if $j<s$) 
or $d(z_{m_j}, w_{m_j})\le \frac1{k_{\min\nolimits_j}-1}$ 
(if $j\ge s$); here $k_{\min\nolimits_j}$ is the least number 
among $k_{m_j}, \dots, k_{m_{j+1}-1}$. 
We set 
$$%
\begin{aligned}
&f\restrict{\{2, \dots, m_1-1\}}=f_1, 
&&f(m_1)=k_{\min\nolimits_0},\\ 
&f\restrict{\{m_1+1, \dots, m_2-1\}}=f_2, 
&&f(m_2)=k_{\min\nolimits_1},\\ 
&\dots,\\ 
&f\restrict{\{m_{s-1}+1, \dots, m_{s}-1\}}=f_s, 
&&f(m_s)=k_{\min\nolimits_{s-1}},\\ 
&f\restrict{\{m_s+1, \dots, m_{s+1}-1\}}=f_s, 
&&f(m_{s+1})=k_{\min\nolimits_{s+1}},\\
&f\restrict{\{m_{s+1}+1, \dots, m_{s+2}-1\}}=f_{s+1},\qquad 
&&f(m_{s+2})=k_{\min\nolimits_{s+2}},\\
&\dots,\\ 
&f\restrict{\{m_{r-1}+1, \dots, m_r-1\}}=f_r, 
&&f(m_r)=k_{\min\nolimits_r},\\  
&f\restrict{\{m_{r}+1, \dots, m_{r+1}-1\}}=f_r.
\end{aligned}
$$%
For $i\le 
n-1$, we have 
$(w_i, z_{i+1})\in \mathcal U_{k_i}$ (by assumption), $(z_i, w_i)\in 
\mathcal U_{f(i)-1}$ (by construction), the $k_i$ are different numbers 
larger than~1, the $f(i)$ are different numbers of the form 
$k_j$, $f(i)>k_m$ for all $i$, and $k_i>k_m$ for $i\ne m$.  
By Lemma~\ref{lemma0},  
\begin{multline*}
(w_1, w_m)=(w_1, z_2)\circ(z_2, w_2)\circ (w_2, z_3)\circ \dots 
\circ (w_{m-1}, z_m)\circ (z_m, w_m)\\
\in \mathcal U_{k_1}\circ 
 \mathcal U_{f(2)}\circ U_{k_2}\circ 
\dots\circ U_{k_{m-1}}\circ \mathcal U_{f(m)}\subset 
\mathcal U_{k_m}
\end{multline*}
and, similarly,  
$$%
(z_{m+1}, z_n)\in \mathcal U_{k_m};
$$%
moreover, by assumption, we have 
$$%
(z_{m+1}, w_m) \in \mathcal U_{k_m}.
$$%
Therefore, 
$$%
(w_1, z_n)\in 3\mathcal U_{k_m}=3\mathcal V_{k_m^2}\subset 
\mathcal V_{k_m^2-1}\subset \mathcal V_{(k_m-1)^2}=
\mathcal U_{k_m-1}.
$$%
The word $(z_1-w_1)+ (z_n-w_n)$ satisfies condition $(*)$, 
because the word $\sum_{i=1}^n(z_i-w_i)$ 
satisfies this condition by assumption (see\ Remark~\ref{remark1}); 
hence $d(z_n, w_n)\le \frac1{k_m-1}$ (if 
($6_{\text{left}}$) holds) or $d(z_1, w_1)\le \frac1{k_m-1}$ 
(if ($6_{\text{right}}$) holds).
\end{proof}

\begin{cor}\label{corollary3}
Suppose that, for $j\le k$, 
$\sum_{i=1}^{n_j}(z^{(j)}_i-w^{(j)}_i)$ 
and $\sum_{i=1}^{n_j-1}(z^{(j)}_{i+1}-w^{(j)}_i)$ 
are words satisfying condition $(*)$, 
$x_j=z^{(j)}_1$, and $y_j=w^{(j)}_{n_j}$. 
Suppose also that $\mathcal U_N\circ (x_r, y_s)\circ \mathcal U_N 
\subset \mathcal U_{\frac 1{d(x_r, y_s)}}$ and 
$\mathcal U_N\circ (x_r, y_s)\circ \mathcal U_N 
\cap  \mathcal U_{\frac 1{d(x_r, y_s)}+1}
=\varnothing$ for all $r, s\le k$. 
Finally, suppose that 
$(z^{(j)}_{i+1}, w^{(j)}_i)\in \mathcal U_{k^{(j)}_i}$ 
for any $j\le k$ and $i\le n_i-1$, where  
the $k^{(j)}_i$ are different positive integers larger than~$N+2$ 
for each $j$.  
Then, for any $j\le k$, there exists an $n_0^{(j)}\le n_j$ such that 
$d(x_r, y_s)=d(z_{n_0^{(r)}}, w_{n_0^{(s)}})$ for all $r, s\le k$.  
\end{cor}

\begin{proof}
Take $j\le k$ and consider the word 
$\sum_{i=1}^{n_j}(z^{(j)}_i-w^{(j)}_i)$; for convenience, 
we omit the index $j$. 

Take some $n_0\le n$ for which $d(z_{n_0}, w_{n_0})$ 
is maximal among all $d(z_i, w_i)$ with $1\le i\le n$. 
Suppose for definiteness that $n_0<n$; if  $n_0>n$, then the left-to-right 
argument described below should be replaced by a similar right-to-left 
argument. 
Let $n_1> n_0$ be the minimum number for which 
$d(z_{n_1}, w_{n_1})$ is largest among all $d(z_i, w_i)$ with 
$i=n_0+1, \dots, n$, and let $m_1$ be such that $k_{m_1}$ is minimal among 
all $k_i$ with $i=n_0, \dots, n_1-1$. Applying Lemma~\ref{lemma3} to 
the word 
$$%
z_{n_0}-w_{n_0}+z_{n_0+1}-w_{n_0+1}+\dots+z_{n_1}-w_{n_1},
$$%
we obtain a one-to-one map 
$$%
f_1\colon \{n_0+1, \dots, n_1-1\}\to \{k_{n_0}, \dots, k_{n_1-1}\}\setminus 
\{k_{m_1}\}
$$%
and the inequalities 
$$%
d(z_{n_1}, w_{n_1})\le \frac1{k_{m_1}-1}\qquad\text{and}\qquad 
d(z_{n_i}, w_{n_i})\le \frac1{f_1(i)-1}
$$%
for $i=n_0+1, \dots, n_1-1$. 
Moreover, by assumption, $d(z_{i+1}, w_i)\le \frac 1{k_i}$ 
for $i=n_0, \dots, n_1-1$.  Therefore, by Lemma~\ref{lemma0},  
\begin{multline*}
(w_{n_0}, w_{n_1})
=
(w_{n_0}, z_{n_0+1})\circ (z_{n_0+1}, w_{n_0+1})\circ \dots 
\circ 
(z_{n_1-1}, w_{n_1-1})\circ (w_{n_1-1}, z_{n_1})\circ(z_{n_1}, 
w_{n_1})\\
\in
\mathcal U_{k_{n_0}}\circ \mathcal U_{f_1(n_0+1)-1}\circ \dots 
\circ \mathcal U_{f_1(n_1-1)-1}\circ \mathcal U_{k_{n_1-1}}\circ 
\mathcal U_{k_{m_1}-1}
\subset \mathcal  U_{k_{m_1}-2}.
\end{multline*}

Consider the word 
$$%
z_{n_1}-w_{n_1}+z_{n_1+1}-w_{n_1+1}+\dots+z_{n_2}-w_{n_2},
$$%
where $n_2> n_1$ is the least number for which 
$d(z_{n_2}, w_{n_2})$ is maximal among all $d(z_i, w_i)$ with 
$i=n_1+1, \dots, n$, and let $m_2$ be such that $k_{m_2}$ is minimal among 
all $k_i$ with $i=n_1, \dots, n_2-1$. Arguing as above, we obtain 
$$%
(w_{n_1}, w_{n_2})\in \mathcal U_{k_{m_2}-2}.
$$%
In the end, 
we join the letters $w_{n_0}$ and $w_n$ by a chain 
$$%
(w_{n_0}, w_n)=
(w_{n_0}, w_{n_1})\circ (w_{n_1}, w_{n_2})\circ \dots \circ 
(w_{n_{t-1}}, w_{n_t})\circ (w_{n_t}, w_{n}),
$$%
where $(w_{n_i-1}, w_{n_i})\in \mathcal U_{k_{m_i}-2}$ for 
$i=1, \dots, t$ and all numbers $m_i$ (and, therefore, 
$k_{m_i}$) are different.  By assumption, $k_{m_i}>N+2$ and 
$w_n=y$; hence Lemma~\ref{lemma0} implies 
$$%
(w_{n_0}, y)\in \mathcal U_{N}.
$$%
Similarly, 
$$%
(x, z_{n_0})\in \mathcal U_{N}.
$$%

Thus, we have shown that, 
for each $j\le k$, there exists an $n_0^{(j)}\le n_j$ such that 
$(x_j, z_{n^{(j)}_0})\in \mathcal U_{N}$ 
and 
$(w_{n^{(j)}_0}, y_j)\in \mathcal U_{N}$. 
This means that 
$$%
(z_{n^{(r)}_0}, w_{n^{(s)}_0}) \in \mathcal U_N\circ (x_r, y_s) 
\circ \mathcal U_N
$$%
for any $r, s\le k$, which immediately implies 
the required assertion.
\end{proof}

\begin{rem}\label{addrem}
In Corollary~\ref{corollary3}, if $z_1^{(j)}= w^{(j)}_1$, then 
$n_0^{(j)}\ne 1$,  and if $z_{n_j}^{(j)}= w^{(j)}_{n_j}$, then 
$n_0^{(j)}\ne n_j$.
\end{rem}

\begin{lem}\label{lemma4}
Suppose that $g=\sum_{i=1}^k(a_i-b_i)$ is an irreducible word\textup; 
$h=\sum_{i=1}^l(u_i-v_i)$ is an irreducible word 
satisfying condition $(*)$\textup; 
$\mathcal U_N\circ (a_i, b_j)\circ \mathcal U_N \subset 
\mathcal U_{\frac 1{d(a_i, b_j)}}$ and 
$\mathcal U_N\circ (a_i, b_j)\circ \mathcal U_N \cap 
\mathcal U_{\frac 1{d(a_i, b_j)+1}}=\varnothing$ for $i, j\le k$\textup; 
$(u_i, v_i) \in \mathcal U_{(N+3)\cdot i}$ for $i\le l$\textup; 
a decomposition $g+h=\sum_{i=1}^m(z_i-w_i)$ 
is irreducible and satisfies condition $(*)$\textup; 
and $(z_i, w_i) \in \mathcal U_{n\cdot i}$ for $i\le m$.
Then there exists a decomposition 
$g=\sum_{i=1}^k(x_i-y_i)$ satisfying $(*)$ 
in which $(x_i, y_i) \in \mathcal 
U_{n\cdot i}$ for $i\le k$.
\end{lem}

\begin{proof}
The decomposition $g+h=\sum_{i=1}^m(z_i-w_i)$ is obtained from 
$\sum_{i=1}^k(a_i-b_i) + \sum_{i=1}^l(u_i-v_i)=g+h$ by 
canceling pairs of equal letters with opposite signs. 
We assume that the cancellations are fixed and each letter in 
this decomposition remembers to which word ($g$ or $h$) it belonged 
before cancellation and which position in this word it occupied. 
In other words, when we say, e.g., that $z_i$ is a letter from 
$g$, this does not merely means that $z_i$ equals some letter $a_j$; 
this means also that some letter $a_j$ from the word $g$ has not been 
canceled in $\sum_{i=1}^k(a_i-b_i) + \sum_{i=1}^l(u_i-v_i)$ (while 
some other letter equal to $a_j$ might have been canceled)
and has become the letter $z_i$.  Possibly, some other letter 
$z_r$ also equals $a_j$, but $z_r$ is not $a_j$, 
because $a_j$ is $z_i$; 
this letter $z_r$ is some other letter $a_s$, or even 
a letter from $h$. To emphasize that, considering 
letters of  $g+h=\sum_{i=1}^m(z_i-w_i)$,  we mean 
letters together with their origins, we use the sign $\equiv$ 
instead of $=$; thus, in the above example,  $z_i\equiv a_j$ but 
$z_r\not\equiv a_j$ (although $z_r = a_j$). 

Take any letter $x_1$ included in the word $g$ with coefficient 
1 (e.g.,  $x_1\equiv a_1$). Our immediate goal is to 
define a letter $y_1$. For this purpose, we shall 
construct a chain of letters of the forms $z_i\equiv v_j$ and 
$w_i\equiv u_j$ until we reach a letter
from $g$; this letter will be $y_1$.

\emph{Link 1.} If $x_1$ is not canceled in the word 
$g+h$, then $x_1\equiv z_{i_1}$ for some $i_1\le m$. If the 
corresponding letter $-w_{i_1}$ is a letter from $g$, then we set 
$y_1\equiv w_{i_1}$; otherwise (i.e., if this is a letter from 
$h$), we have $w_{i_1}\equiv v_{j_1}$ for some $j_1\le l$.  If the 
letter $x_1$ is canceled in the word $g+h$, then it is canceled by 
a letter from $h$ (because $g$ is irreducible), i.e., by 
$-v_{j_1}$ for some $j_1\le l$. We have either found 
$y_1$ or defined $v_{j_1}$ and (possibly) $w_{i_1}\equiv 
v_{j_1}$ and $z_{i_1}$. 

\emph{Link 2.} If the letter $u_{j_1}$ corresponding to 
the $v_{j_1}$ found at the preceding step 
is not canceled in the word $g+h$, 
then $u_{j_1}\equiv z_{i_2}$ for some $i_2\le m$. If the 
corresponding $-w_{i_2}$ is a letter from $g$, then we set $y_1 
\equiv w_{i_2}$; otherwise, we have $w_{i_2}\equiv v_{j_2}$ 
for some $j_2\le l$. If the letter $u_{j_1}$ is canceled in 
the word $g+h$, then it is necessarily canceled by a 
letter $b_\alpha$ from $g$, and we take this letter for $y_1$; then 
$y_1\equiv b_\alpha = u_{j_1}$. We have either found $y_1$ or defined 
$z_{i_2}\equiv u_{j_1}$ and $w_{i_2}\equiv v_{j_2}$. 

Continuing, we obtain $y_1$ in the end. 

Applying this procedure to all letters of $g$ with positive coefficients 
in turn, we 
obtain a partitioning of the letters of $g$ into pairs $x_s, y_s$ 
together with chains of letters 
$$%
z'_{i_1(s)} \equiv x_s,\ 
w'_{i_1(s)} \equiv v_{j_1(s)},\ 
z'_{i_2(s)}\equiv u_{j_1(s)},\  
w'_{i_2(s)}\equiv v_{j_2(s)},\  
\dots\,,\ 
z'_{i_{r_s}(s)}\equiv u_{j_{r_s-1}(s)},\ 
w'_{i_{r_s}(s)}\equiv y_s,\ 
$$%
where $z'_{i_1(s)}\equiv z_{i_1(s)}$ (if $x_s$ is not canceled in 
$g+h$; in this case, $w'_{i_1(s)}\equiv w_{i_1(s)}$) 
or $z'_{i_1(s)}= v_{j_1(s)}$ 
(if $x_s$ is canceled by $-v_{j_1(s)}$), 
$z'_{i_t(s)}\equiv z_{i_t(s)}\equiv u_{j_t(s)}$ and 
$w'_{i_t(s)}\equiv w_{i_t(s)}\equiv v_{j_t(s)}$ for $t=2, \dots, r_s-1$, 
$z'_{i_{r_s}(s)}\equiv u_{j_{r_s}(s)}$, 
and $w'_{i_{r_s}(s)}\equiv w_{i_{r_s}(s)}$ or 
$w'_{i_{r_s}(s)}\equiv b_\alpha= u_{j_{r_s-1}(s)}$ for some $\alpha_s$. 
The sets $\{i_\alpha(s)\}$ are disjoint 
for different $s$. 
The sums 
$\sum_{t=1}^{r_s}(z'_{i_t(s)}-w'_{i_t(s)})$ 
satisfy the conditions of Corollary~\ref{corollary3}. 
Indeed, these sums satisfy condition $(*)$, because 
their terms are divided into the pairs $z'_{i_t(s)}-w'_{i_t(s)}$, 
which belong to a decomposition of $g+h$ 
satisfying condition $(*)$.  The first and last pairs may 
differ from the corresponding terms of the decomposition 
of $g+h$, but they equal zero (the empty word) in this case; i.e., 
either  $z'_{i_1(s)} = z_{i_1(s)}$ and $w'_{i_1(s)} = w_{i_1(s)}$ 
or $z'_{i_1(s)}=w'_{i_1(s)}$, and 
either $z'_{i_{r_s}(s)}= z_{i_{r_s}(s)}$ and 
$w'_{i_{r_s}(s)}= w_{i_{r_s}(s)}$ or $w'_{i_{r_s}(s)}= 
z'_{i_{r_s}(s)}$; so, condition $(*)$ is not 
violated. The sums $\sum_{t=1}^{r_s-1}(z'_{i_{t+1}(s)}-w'_{i_t(s)})$ 
also satisfy condition $(*)$, because each pair 
$z'_{i_{t+1}(s)}-w'_{i_t(s)}=u_{j_t(s)}-v_{j_t(s)}$ is contained in 
a decomposition of $h$ satisfying condition $(*)$.  Moreover, 
by assumption, we have 
$(z'_{i_{t+1}(s)},w'_{i_t(s)})= (u_{j_t(s)},v_{j_t(s)})\in 
\mathcal U_{(N+3)\cdot j_t(s)}$, and the $(N+3)\cdot j_t(s)$ 
are different numbers larger than $N+2$. Finally, since 
all $x_r$ and $y_s$ are letters of the word 
$g=\sum_{i=1}^k(a_i-b_i)$, it follows from the conditions of the 
lemma being proved that the remaining condition of 
Corollary~\ref{corollary3} holds too; namely, $\mathcal U_N\circ 
(x_r, y_s)\circ \mathcal U_N \subset \mathcal U_{\frac 1{d(x_r, 
y_s)}}$ and $\mathcal U_N\circ (x_r, y_s)\circ \mathcal U_N \cap 
\mathcal U_{\frac 1{d(x_r, y_s)}+1} =\varnothing$ for all $r, s\le 
k$.  Therefore, for all $s$, there exist $n_0^{(s)}\in \{i_t(s): 
t=1, \dots, r_s\}$  such that 
$d(x_r, y_t)=d(z'_{n_0^{(r)}}, w'_{n_0^{(t)}})$ for any $r, t\le k$, 
and the numbers $n_0(s)$ are different for different $s$ (because the 
sets $\{i_t(s): t=1, \dots, r_s\}$ are disjoint). By Remark~\ref{addrem}, 
$n_0(s)\ne i_1(s)$ if $z'_{i_1(s)} \ne z_{i_1(s)}$ 
or $w'_{i_1(s)} \ne w_{i_1(s)}$  (i.e., $z'_{i_1(s)}=w'_{i_1(s)}$) 
and  $n_0(s)\ne i_{r_s}(s)$ if $z'_{i_{r_s}(s)} \ne z_{i_{r_s}(s)}$ 
or $w'_{i_{r_s}(s)} \ne w_{{r_s}_1(s)}$  
(i.e., $z'_{i_{r_s}(s)}=w'_{i_{r_s}(s)}$). Thus, we have 
$d(x_r, y_t)=d(z_{n_0^{(r)}}, w_{n_0^{(t)}})$ for $r, t\le k$.   
Since the sum $\sum_{i=1}^m(z_i-w_i)$ satisfies condition $(*)$, 
it follows that the sum $\sum_{t=1}^k(z_{n_0^{(t)}}-w_{n_0^{(t)}})$ 
also satisfies condition $(*)$ (see Remark~\ref{remark1}); 
hence $g=\sum_{s=1}^k (x_s-y_s)$ satisfies condition~$(*)$.  
Finally, it follows from  
$d(x_i, y_i) = d(z_{n_0^{(i)}}, w_{n_0^{(i)}})$ that $(x_i, y_i)\in 
\mathcal U_{n_0^{(i)}\cdot n}$ for $i\le k$. Since all $n_0^{(i)}$ are 
different, 
we can assume that each $(x_i, y_i)$ belongs to $\mathcal U_{i\cdot n}$ 
(otherwise, 
we renumber the terms $x_i -y_i$ and recall that $\mathcal U _r\supset 
\mathcal U_s $ for  $r\le s$). 
\end{proof}

All is ready for  the proof of the last assertion, from which 
Theorem~1 follows immediately. 

Recall that, at the beginning of the paper, we defined the sets 
$W_n(\boldsymbol{\mathcal U})$, which form a neighborhood base at 
zero for a (metrizable) group topology $\mathcal T_{\boldsymbol 
{\mathcal U}}$ on $A(X)$. We set 
\begin{multline*}
W^*_n(\boldsymbol{\mathcal U}) = 
\bigcup_{k\in \omega} \bigl\{\sum_{i=1}^k(x_i-y_i):
(x_i, y_i)\in \mathcal U_{n\cdot i}, \\
\text{the decomposition $\sum_{i=1}^k(x_i-y_i)$ 
satisfies condition $(*)$}\bigr\}. 
\end{multline*}

\begin{claim}
\begin{enumerate}
\item[\upshape (i)] $W^*_{n}(\boldsymbol {\mathcal U})\subset 
W_n(\boldsymbol {\mathcal U})$ for all $n$\textup;
\item[\upshape (ii)] $W_{n+1}(\boldsymbol {\mathcal U})\subset 
W^*_n(\boldsymbol {\mathcal U})$ for all $n$\textup;
\item[\upshape (iii)] for any $n\in \omega$ and any 
$g\in W^*_{n}(\boldsymbol {\mathcal U})$, 
there exists an $n_0\in \omega$ for which $g+ W_{n_0}(\boldsymbol 
{\mathcal U})\subset W^*_{n}(\boldsymbol {\mathcal U})$\textup; 
\item[\upshape (iv)] for any $k,n\in \omega$ and any word 
$g=\sum_{i=1}^k(a_i-b_i)$, 
there exists an $n_0\in \omega$ such that the condition $g+ 
W^*_{n_0}(\boldsymbol {\mathcal U})\cap W^*_{n}(\boldsymbol 
{\mathcal U})\ne \varnothing$ implies $g\in W^*_{n}(\boldsymbol 
{\mathcal U})$.
\end{enumerate} 
\end{claim}

\begin{proof}
Assertion (i) is obvious; (ii) is Corollary~\ref{corollary1}.  
Assertion~(iii) follows from 
Corollary~\ref{corollary2}. Indeed, suppose that $g=\sum_{i=1}^k (x_i-y_i)$, 
the decomposition $\sum_{i=1}^k (x_i-y_i)$ satisfies condition $(*)$, 
and $(x_i, y_i)\in \mathcal U_{n\cdot i}$ for $i\le k$. We can 
assume that $x_i\ne y_i$ for $i\le k$, because if $x_j=y_j$ for some 
$j$, then we can delete the term $x_j-y_j$ from the sum 
$\sum_{i=1}^k (x_i-y_i)$; i.e., we can set $x'_i=x_i$ and $y'_i= 
y_i$ for $i< j$ and $x'_i=x_{i+1}$ and $y'_i=y_{i+1}$ for $i=j, 
\dots, k-1$; we have $g=\sum_{i=1}^{k-1} (x'_i-y'_i)$, 
the decomposition $\sum_{i=1}^{k-1} (x'_i-y'_i)$ satisfies condition 
$(*)$ (see Remark~\ref{remark1}), and $(x'_i, y'_i)\in \mathcal 
U_{n\cdot i}$ for $i\le k-1$ (for $i\ge j$, we have $(x'_i, y'_i)\in 
\mathcal U_{n\cdot {i+1}}\subset \mathcal U_{n\cdot i}$).  Thus, suppose 
that $x_i\ne y_i$; in this case, 
the decomposition $\sum_{i=1}^k (x_i-y_i)$ is irreducible 
(i.e., $x_i\ne y_j$ for any $i, j\le k$),
because it satisfies $(*)$.
Since all $\mathcal U_r$ are 
clopen and form a base for a uniformity generating the initial 
(completely regular) topology on $X$, we 
can find $N$ for which the conditions  
of Corollary~\ref{corollary2}  hold; after that, it remains to 
set $n_0=N+2$:  if $h\in W_{N+2}(\boldsymbol{\mathcal U})$, then 
$h\in W_{N+1}(\boldsymbol{\mathcal U})$ (see (ii)) and, by 
Corollary~\ref{corollary2}, $g+h\in W^*_n(\boldsymbol{\mathcal U})$.  
Assertion~(iv) is derived from Lemma~\ref{lemma4} in a  similar 
way ($n_0=N+4$).  
\end{proof}

It follows from (i)--(iii) that the sets $W^*_n(\boldsymbol{\mathcal 
U})$ are open in the topology $\mathcal T_{\boldsymbol{\mathcal 
U}}$ and form a neighborhood base at zero for this topology; 
(iv) says that each $W^*_n(\boldsymbol{\mathcal U})$ is closed in 
$\mathcal T_{\boldsymbol{\mathcal U}}$.

\begin{rem}\label{rem-graev}
Let $\rho$ be a metric on $X$  such that 
$\mathcal U_{i}\subset \{(x,y):\rho(x,y)\le 
\frac1{2^{i}}\}\subset \mathcal U_{i-1}$ for any $i\ge 1$ (it 
exists by Theorem~8.1.10 from~\cite{Eng}). Then the topology on 
$A(X)$ generated  by the Graev extension of $\rho$ is no 
stronger than $\mathcal T_{\boldsymbol{\mathcal U}}$. 
Indeed, if $g\in 
W_n(\boldsymbol{\mathcal U})$, then $g= \sum_{i=1}^k(x_i-y_i)$, where   
$(x_i, y_i)\in \mathcal U_{n\cdot i}$ for $i\le k$, and 
$\sum_{i=1}^k\rho(x_i-y_i)\le \sum_{i=1}^k \frac1{2^{n\cdot i}} 
< \frac 1{2^n}$. Since the Graev norm $\|g\|_\rho$ of the element $g$ 
is defined as $\min\Bigl\{\sum_{i=1}^m\rho(u_i,v_i):\ 
m\ge 1, \ g= \sum_{i=1}^m(u_i-v_i)\Bigr\}$, we have 
$\|g\|_\rho < \frac 1{2^n}$.  Thus, each Graev ball of radius $\frac 
1{2^n}$ centered at zero contains some base neighborhood 
$W_n(\boldsymbol{\mathcal U})$ of zero in the topology  $\mathcal 
T_{\boldsymbol{\mathcal U}}$.  Since the space $X$ is closed 
in the free group with the Graev topology, it is also closed 
in the free group with the topology $\mathcal 
T_{\boldsymbol{\mathcal U}}$. 
\end{rem}

\section{A Metrizable Group 
with Noncoinciding Dimensions}

We denote the Cantor set $2^\omega$ 
by $C$. 
The elements of $C$ are infinite 
sequences of zeros and ones. The topology of $C$ has 
a standard base, which is a tree under 
inclusion; the $n$th-level elements of this tree are sets 
of sequences whose first $n$ members coincide; different 
elements of the same level do not intersect. Clearly, all base
neighborhoods of the same point of $C$ are comparable, and 
larger neighborhoods belong to levels with smaller numbers.  
We denote the elements of the Cantor set $C$  itself by the letters $x$, 
$y$, $z$, \dots and the infinite sequences of such elements  (i.e., 
the elements of the set $C^\omega$) by the same letters in 
boldface:  $\xxx$, $\yyy$, $\zzz$, \dots; we denote the value of  
a sequence $\xxx $ at $n$ by $\xxx (n)$. The 
restriction of a sequence $x\in C$ to $\{0, \dots, n-1\}$ (i.e., the ordered 
set of the first $n$ elements of this sequence) is denoted by 
$x |_n$. Thus, the $n$th-level elements of the base-tree have the form 
$\{y\in C: y|_{n+1}= x|_{n+1}\}$ for $x\in C$.

By $I$ we denote the usual interval $[0, 1]$. Let $t\in I$. 
If $t=(2k+1)/2^{n}$ for some 
positive integers $k$ and $n$, then we define the order $t$ as $\ord 
t = n$.  We assume that $\ord 0 =\ord 1 = 0$. For all other 
numbers $t\in [0, 1]$, we set $\ord t=\infty$. 

For $n\in \omega$, we define the neighborhood $I_n(t)$ of a 
number $t\in (0, 1)$ to be the interval $I_n(t)=(a_n(t), b_n(t))$, 
where $a_n(t)$ and $b_n(t)$ are the dyadic rationals of 
minimal order for which $b_n(t)-a_n(t)=1/2^{n}$ and $t\in (a_n(t), 
b_n(t))$; we set $I_n(0)=[0, 1/2^{n+1})$ and $I_n(1)=(1-1/2^{n+1}, 
1]$. Thus, if $0<\ord t\le n$, i.e., $t=k/2^n$ for 
some (possibly, even) $k$, then $a_n(t)=(2k-1)/2^{n+1}$ and 
$b_n(t)=(2k+1)/2^{n+1}$ (and hence $\ord a_n(t)=\ord b_n(t)=n+1$), 
and if $\ord t>n$, then  $a_n(t)=k/2^{n}$ and $b_n(t)=(k+1)/2^{n}$ 
for some $k$ (and hence the order of one of the numbers $a_n(t)$ and 
$b_n(t)$ equals $n$ and the order of the other is strictly less 
than $n$). 

Let $A\subset C$. We set 
$$%
\nu\mu_0(A)=\{(\xxx , t)\in C^\omega\times I:
\text{ $\xxx (n) \in A$ for $n\ne \ord t$,} 
\text{ $\xxx (n) \in C\setminus A$ for $n=\ord t$}\}
$$%
and endow $\nu\mu_0(A)$ with the topology generated by the sets 
of the form 
$$
U_n(\xxx , t)=
\begin{cases}
\{(\yyy , s)\in \nu\mu_0(A):
\text{$s\in I_n(t)$, $\yyy (i)=\xxx (i)$ 
for $i\le n$}\}&
\text{if $n<\ord t$},\\
\{(\yyy , s)\in \nu\mu_0(A):
\text {$s\in I_n(t)$; 
$\yyy (i)=\xxx (i)$ 
for $i\le n+1$, $i\ne \ord t$;}&\\
\qquad\qquad\qquad\qquad \text{$\yyy (i)|_{n+1} = \xxx (i)|_{n+1}$
for $i=\ord t$}\}&
\text{if $\ord t\le n$}.  
\end{cases}  
$$
According to Mrowka~\cite{1}, the space $\nu\mu_0(A)$ is metrizable 
and 
$\ind\nu\mu_0(A) =0$; moreover, if $A$ is everywhere 
dense in $C$ and the set $C\setminus A$ is of second 
category, then $\dim\nu\mu_0(A) >0$. 

The projection $\pi (\nu\mu_0(A))$ of the set $\nu\mu_0(A) 
\subset C^\omega\times I$ on the first factor consists of 
all sequences $\xxx \in C^\omega$ each of which takes at most one 
value not in $A$.

For $A$ we take the set $\sigma 2^\omega$ of binary 
sequences with only finitely many elements different 
from 0. For each nonzero $x\in A$,  we define its length $\len x$ to 
be the number of the last nonzero term of the sequence $x$; 
we set $\len 0=0$ (thus, $\len 0100\dots = 2$). 

For $\xxx \in\pi (\nu\mu_0(A))$ and $n,i\in \omega$, we fix 
a maximal base neighborhood $\J_{n}(\xxx, i)$ of $\xxx (i)$ 
of level $\ge n$ such that 
\begin{enumerate} 
\item if $\xxx (j) \in A$ for all 
$j\le n$, then the lengths of all elements of the intersection 
$\J_{n}(\xxx,i)\cap A$ (except, possibly, 
the point $\xxx (i)$ itself) are larger than all lengths $\len \xxx 
(j)$ for $j\le n$; 
\item if $\xxx (j) \notin A$ for some $j\le n$, 
then the lengths of all elements of the intersection 
$\J_{n}(\xxx, i)\cap A$ (except, possibly, the point 
$\xxx(i)$ itself) are larger than the lengths $\len \xxx (j)$ 
for all $j\le n+1$ such 
that $\xxx (j)\in A$.  
\end{enumerate} 
Since all sets  
of the form $\J_{n}(\xxx, i)$ are elements of the base-tree, 
it follows that, for any $\xxx ,\yyy \in\pi (\nu\mu_0(A))$ and 
any $n,k,i,j\in \omega$, either the sets $\J_{n}(\xxx, i)$ and $\J_{k}(\yyy, 
j)$ are disjoint or one of them is contained in 
the other.

For $(\xxx , t)\in \nu\mu_0(A)$, we set 
$$%
V_n(\xxx , t)=
\begin{cases}
\{(\yyy , s)\in \nu\mu_0(A):
\text{$s\in I_n(t)$, $\yyy (i)=\xxx (i)$ 
for $i\le n$, $\yyy (n+1)\in \J_{n}(\xxx,n+1)$}\}&
\text{if $n<\ord t$},\\ 
\{(\yyy , s)\in \nu\mu_0(A):  \text{$s\in I_n(t)$; 
$\yyy (i)=\xxx (i)$ 
for $i\le n+1$, $i\ne \ord t$;}&\\
\qquad\qquad\qquad\qquad\text{$\yyy (i)\in \J_{n}(\xxx, i)$
for $i=\ord t$}\}&
\text{if 
$\ord t\le n$}.  
\end{cases} 
$$%
Clearly, the sets of the form $V_n(\xxx , t)$ constitute a base 
for the topology of $\nu\mu_0(A)$. 

\begin{rem}\label{mrorem}
Suppose that $(\xxx , t), (\yyy , s)\in \nu\mu_0(A)$, $n,n'\in 
\omega$, and $I_{n'}(s)\cap \{r\in [0, 1]:\ord r\le n+1,\ r\ne 
s\}=\varnothing$ (this implies, in particular, that $n'\ge n$). Then 
one of following four cases occurs:  
\begin{enumerate} 
\item[\upshape (i)] 
$(\yyy , s)\in V_n(\xxx , t)$; 
\item[\upshape (ii)] 
$V_{n'}(\yyy , s)\cap V_{n}(\xxx , t)=\varnothing$; 
\item[\upshape (iii)] 
$s\in I_n(t)$, $\ord s, \ord t>n$, 
$\yyy (i)= \xxx (i)$ for all $i\le n$, and $\xxx (n+1)\in \J_{n'}(\yyy, 
n+1)$; moreover, in this case, 
$V_{n}(\xxx , t)\subset V_{n}(\yyy , s)$; 
\item[\upshape (iv)] 
$s\in \overline {I_n(t)}\setminus I_n(t) = 
\{a_n(t), b_n(t)\}$ and $\yyy (i)=\xxx (i)$ for all $i\le n$ such 
that $i\ne \ord t, \ord s$.  
\end{enumerate}

Indeed, if $s\notin \overline {I_n(t)}$, then 
$I_{n'}(s)\cap I_{n}(t)=\varnothing$ and 
$V_{n'}(\yyy , s)\cap V_{n}(\xxx , t)=\varnothing$, i.e., 
condition (ii) holds.

If $s\in I_n(t)$, $\ord s\le n$ (this can happen only if 
$s=t$), and $\yyy (i)\ne \xxx (i)$ for some $i\le n+1$ such that 
$i\ne \ord t$, then (ii) holds.

If $s\in I_n(t)$, $\ord s=k\le n$ (then $s=t$),  
and $\yyy (i)= \xxx (i)$ for all $i\le n+1$ such that 
$i\ne k$, then either (a)~$\yyy (k)\in \J_n
(\xxx,k)$ (and then (i) holds), (b)~$\J_{n'}(\yyy,k)
\cap \J_n(\xxx, k)=\varnothing$ (and then  (ii) holds), 
or (c)~$\J_{n'}(\yyy, k)\supset \J_{n}(\xxx, k)$ and 
$\yyy (k)\notin \J_n(\xxx,k)$. In case~(c), 
$\J_{n'}(\yyy, k)$ is 
a base neighborhood of the point $\xxx (k)$, its 
level is at least $n'\ge n$, and the lengths of all 
elements of the intersection $\J_{n'}(\yyy, k)\cap A$ 
are larger than the length $\len \xxx(j)=\len \yyy (j)$ 
for all $j\le n+1$ such that $\xxx(j)=\yyy (j)\in A$ 
(the points $\xxx (k)$ and $\yyy(k)$ themselves do not 
belong to $A$). This contradicts the maximality of the neighborhood 
$\J_n(\xxx, k)$. 

Suppose that $s\in I_n(t)$, $\ord s > n$, and $\ord t=k\le n$. Then 
$\yyy(k) \in A$, $\xxx (k)\notin A$, and $\xxx (n+1)\in A$. If 
(ii) does not hold, then there exists a $(\zzz, h)\in 
V_{n'}(\yyy , s)\cap V_{n}(\xxx , t)$. By the definition of 
$V_{n'}(\yyy , s)$ and $V_{n}(\xxx , t)$, we have $\zzz(i)= 
\yyy (i)=\xxx (i)$ for all $i\le n$ 
different from $k$, $\zzz(k)=\yyy (k)\in \J_n(\xxx, k)$, and 
$\zzz(n+1)=\xxx (n+1)\in \J_{n'} (\yyy,n+1)$.  
We have $k\le n$, $\xxx(k)\notin A$, $\yyy(k)\ne \xxx(k)$ (because 
$\yyy(k)\in A$), 
and $\xxx(n+1)\in A$; thus, it follows from 
$\yyy (k)\in \J_n(\xxx,k)$ that  
$\len \yyy (k)>\len \xxx (i)$ for all $i\le n+1$ different from 
$k$ (in particular, $\len \yyy(k)> \len \xxx (n+1)$). The 
inclusion $\zzz(n+1)=\xxx (n+1)\in \J_{n'} (\yyy,n+1)$ 
implies that either $\len \xxx(n+1)>\len 
\yyy (i)$ for all $i\le n$ or $\xxx (n+1)= \yyy (n+1)$.  The 
former inequality cannot hold, because $\len \xxx (n+1)< \len 
\yyy (k)$; hence $\xxx (n+1)= \yyy(n+1)$.  Thus, $\yyy 
(i)=\xxx (i)$ for $i\le n+1$, $i\ne k$, and $\yyy (k)\in \J_n(\xxx, 
k)$.  This means that (i) holds. 

If $s\in I_n(t)$, $\ord s > n$, $\ord t>n$, 
and $\yyy (i)\ne \xxx (i)$ for some $i\le n$, 
then (ii) holds.

If $s\in I_n(t)$, $\ord s > n$, $\ord t> n$  
and $\yyy (i)= \xxx (i)$ for all $i\le n$, 
then either (a)~$\yyy (n+1)\in \J_n
(\xxx,n+1)$ (and hence (i) holds), (b)~$\J_{n'}(\yyy,n+1)
\cap \J_n(\xxx, n+1)=\varnothing$ (then  (ii) holds), 
or (c)~$\J_{n'}(\yyy, n+1)\supset \J_{n}(\xxx, n+1)$ 
(i.e., (iii) holds). The inclusion $V_{n}(\xxx , t)\subset 
V_{n}(\yyy , s)$ follows from the obvious inclusion 
$\J_{n'}(\yyy , n+1)\subset \J_n(\yyy , n+1)$ (which is an 
immediate consequence of $n'\ge n$).

If 
$s\in \overline {I_n(t)}\setminus I_n(t)$ 
and $\yyy (i)\ne\xxx (i)$ for some $i\le n$ such that 
$i\ne \ord t, \ord s$, then (ii) holds. 
\end{rem}

\begin{claim}\label{emptyboundary}
For any $n\in \omega$, the set
$$%
\mathcal U_n=\bigcup 
\{V_n(\xxx , t) \times V_n(\xxx , t):
(\xxx , t) \in \nu\mu_0(A)\}
$$%
has empty boundary.
\end{claim}

\begin{proof}
Suppose that $(\yyy , s), (\zzz , r)\in \nu\mu_0(A)$, and $n\in 
\omega$.  Take $n'\in \omega$ such that 
$$%
I_{n'}(s)\cap \{t\in [0, 1]:\ord t\le n+1, t\ne s\}=\varnothing,
$$%
$$%
I_{n'}(r)\cap \{t\in [0, 1]:\ord t\le n+1, t\ne r\}=\varnothing,
$$%
and if $\yyy (i)\ne \zzz (i)$ for $i\le n+1$, then 
$$%
\J_{n'}(\yyy,i)\cap \J_{n'}(\zzz,i)=\varnothing.
$$%
Suppose that $V_{n'}(\yyy , s)\times V_{n'}(\zzz , r) \cap 
\mathcal U_n\ne \varnothing$ but 
$((\yyy , s), (\zzz , r))\notin \mathcal U_n$. This means 
that there exist $(\xxx , t)\in \nu\mu_0(A)$, 
$(\yyy ', s')\in V_{n'}(\yyy , s)$, and 
$(\zzz ', r')\in V_{n'}(\zzz , r)$ such that 
$(\yyy ', s')\in V_{n}(\xxx , t)$, 
$(\zzz ', r')\in V_{n}(\xxx , t)$, and either 
$(\yyy , s)\notin V_{n}(\xxx , t)$ or 
$(\zzz , r)\notin V_{n}(\xxx , t)$. For definiteness, 
suppose that 
$(\yyy , s)\notin V_{n}(\xxx , t)$. Then 
(iii) or (iv) from Remark~\ref{mrorem} holds. Suppose that (iv) holds. 
There are the following possibilities: 
\begin{enumerate}
\item 
$\ord t=k\le n$. In this case, $\ord s=n+1$ and $\yyy (i)=
\yyy'(i)\in A$ for $i\le n$. Moreover, 
$\yyy '(n+1)=\xxx (n+1)\in A$ and $\yyy '(k)\in 
\J_n(\xxx,k)$ (because $(\yyy ', s')\in V_n(\xxx , t)$). 
Therefore, $\len \yyy '(k) > \len \xxx (j)$ for all $j\le n+1$ 
different from $k$ (in particular, 
$\len \yyy '(k) > \len \xxx (n+1)$). On the other hand,  $\yyy 
'(n+1)\in \J_{n'}(\yyy,n+1)$ (because $(\yyy ', s')\in V_{n'} (\yyy, 
s)$) and $\yyy '(n+1)\ne \yyy (n+1)$ (because $\ord s=n+1$ and, 
therefore, $\yyy (n+1)\notin A$).  Hence $\len \yyy'(n+1) = \len \xxx 
(n+1)>\len \yyy (k) = \len\yyy '(k)$.  This is 
impossible. 
\item $\ord t> n$. In this case, 
$\ord s = k\le n$.  Suppose that $r\ne s$.  
If $r\notin \overline I_n(t)$, 
then $I_{n'}(r)\cap I_{n}(t)=\varnothing$ and  
$V_{n'}(\zzz , r)\cap V_n(\xxx , t)=\varnothing$, 
which contradicts the assumption.  
Therefore, $r\in \overline I_n(t)$, and $\ord r\ne \ord s=k$ (the 
endpoints of the interval $I_n(t)$ are of different orders, and all 
interior points of this interval have orders larger than $n$).  
Thus, $\yyy (k)\in C\setminus A$, whereas $\zzz 
(k)\in A$.  The number $n'$ was chosen so that 
$\J_{n'}(\yyy, k)\cap \J_{n'}(\zzz,k) =\varnothing$; in 
particular, $\zzz (k)\notin \J_{n'}(\yyy,k)$.  Since $(\yyy ', 
s')\in V_{n'}(\yyy , s)$, $(\zzz ', r')\in V_{n'}(\zzz , r)$, and, 
moreover, $\ord s=k$, $\ord r\ne k$, and $k\le n$, it follows 
that $\yyy '(k)\in \J_{n'}(\yyy,k)$ and $\zzz '(k) =\zzz (k)\notin 
\J_{n'}(\yyy,k)$.  Therefore, $\yyy '(k)\ne \zzz '(k)$, and at least 
one of these numbers is not equal to $\xxx (k)$, i.e., at least one 
of the pairs $(\yyy ', s')$ and $(\zzz ', r')$ does not belong 
to the set $V_n(\xxx , t)$, which contradicts the definition of these 
pairs.  Hence $r=s$. The same argument shows that $\yyy (i)=\zzz 
(i)$ for all $i\le n$: if $y(i)\ne z(i)$, then at least one of the 
numbers $y'(i)$ and $z'(i)$ is not equal to $x(i)$, and the 
corresponding pair does not belong to $V_n(\xxx , t)$. 
Since $(\yyy', s')\in V_n(\xxx, t)$ and $\ord t>n$, 
we have $\yyy'(i) = \xxx(i)$ for all $i\le n$; since 
$(\yyy', s')\in V_{n'}(\yyy, s)$ and $\ord s = k\le n$, we have 
$\yyy(k)\notin A$ and $\yyy'(k)=\xxx(k)\in J_{n'}(\yyy, k)$. 
Therefore, $\len\xxx (k)>\len \yyy (i)$ for 
all $i\le n+1$ different from $k$ (in particular, $\len\xxx 
(k)>\len \yyy (n+1)$).  Since $(\yyy ', s') \in V_n(\xxx , t)$, we 
have $\yyy '(n+1)\in \J_n(\xxx, n+1)$. Therefore, either $\len \yyy 
'(n+1)>\len \xxx(k)$ or  $\yyy '(n+1)=\xxx (n+1)$. On the 
other hand, $(\yyy ', s') \in V_{n'}(\yyy , s)$ and $\ord s=k\le n$, 
whence $\yyy '(n+1)=\yyy (n+1)$. Thus, the inequality $\len \yyy'(n+1)>\len 
\xxx (k)$ cannot hold; hence $\yyy(n+1)=\yyy '(n+1)=\xxx 
(n+1)$.  Similarly, $\zzz(n+1)=\xxx (n+1)$.  Thus, $s=r$ and 
$\yyy (i)=\zzz (i)$ for $i\le n+1$; therefore,  $(\zzz , r)\in 
V_n(\yyy , s)$, i.e., $((\yyy , s), (\zzz , r))\in \mathcal U_n$.  
\end{enumerate}

Now, suppose that condition (iii) from Remark~\ref{mrorem} holds. 
If $(\zzz , r)\in V_n(\xxx , t)$, then 
$(\zzz , r)\in V_n(\yyy , s)$ and 
$((\yyy , s), (\zzz , r))\in \mathcal U_n$. Suppose that 
$(\zzz , r)\notin V_n(\xxx , t)$. Since 
$V_{n'}(\zzz , r)\cap V_n(\xxx , t)\ne \varnothing$, 
it follows that one of conditions (iii) and (iv) with $\zzz $ 
instead of $\yyy$ and $r$ instead of $s$ holds. The case in 
which (iv) holds has just been considered. Suppose 
that (iii) holds.  We have  $s,r\in I_n(t)$; $\ord s, \ord 
r, \ord t> n$; $\yyy (i)= \zzz (i)=\xxx (i)$ for all $i\le 
n$ (because $(\xxx, t)\in V_n(\xxx, t)$, $V_n(\xxx, t)\subset 
V_n(\yyy, s)$ by condition (iii) for $(\xxx, t)$ and 
$(\yyy, s)$, and $V_n(\xxx, t)\subset V_n(\zzz, r)$ 
by condition (iii)  for $(\xxx, t)$ and $(\zzz, r)$); 
$\xxx (n+1)\in \J_{n'}(\yyy,n+1)$; and $\xxx (n+1)\in 
\J_{n'}(\zzz,n+1)$.  Therefore, $\J_{n'}(\yyy,n+1)\subset 
\J_{n'}(\zzz,n+1)$ or $\J_{n'}(\zzz,n+1)\subset \J_{n'}(\yyy, 
n+1)$.  For definiteness, suppose that $\J_{n'}(\yyy,n+1)\subset 
\J_{n'}(\zzz,n+1)$.  Then $\yyy (n+1)\in \J_{n'}(\zzz,n+1)
\subset \J_{n}(\zzz,n+1)$. It remains to note that 
$I_n(r)=I_n(t)$ (because $\ord t>n$ 
and $r\in I_n(t)$). This immediately implies $s\in 
I_n(r)$ and $(\yyy , s)\in V_n(\zzz , r)$, i.e., $((\yyy , 
s), (\zzz , r))\in \mathcal U_n$.  This contradiction 
completes the proof.
\end{proof}

It follows immediately from Claim~\ref{emptyboundary} and 
Theorem~\ref{theorem1.1} that the space 
$\nu\mu_0(A)$ can be embedded in a metrizable topological group $G$ 
with $\ind G=0$; moreover, $\nu\mu_0(A)$ is closed in $G$ 
(see~Remark~\ref{rem-graev}). Since $\dim \nu\mu_0(A) > 0$ and 
the group $G$ is metrizable, we have $\dim G>0$. Thus, we have 
obtained an example of a metrizable group with noncoinciding 
dimensions $\ind$ and $\dim$.  

\section[A Zero-Dimensional Metrizable Space]{A Zero-Dimensional 
Metrizable Space Which is not Embedded\\
in a Zero-Dimensional Metrizable Group}

In this section, by a sequence we mean a map 
from an at most countable ordinal to some set and  
consider only sequences with values in $\omega_1$.  
We identify all sequences with ordered sets 
of their values and write them in the form of (finite or infinite) words. 
As in the preceding section, 
we denote sequences by boldface Latin letters,  
but their elements we denote 
by the same letters with subscript-numbers. Thus, the symbol $a_n$ 
always denotes the element number $n$ in the sequence $\aaa$:  
$a_n=\aaa(n)$.  The word whose letters are sequences  
(all but the last must be finite) denotes 
the concatenation of these sequences.  For example, if $\aaa = 
a_0a_1\dots a_n$ and $\bbb = b_0b_1\dots$\,, then $\aaa\bbb = a_0a_1\dots 
a_nb_0b_1\dots$\,.

If $\aaa$ is a sequence of length $\ge n$, then 
$$
\aaa|_n=a_0a_1\dots a_{n-1}
$$
(recall that we assume that 
$\aaa=a_0a_1\dots$); we set $\aaa|_0=\varnothing$. For $m<n$, 
$$
\aaa|^m=a_ma_{m+1}\dots
\qquad\text{and}\qquad
\aaa|_n^m=a_m\dots a_{n-1}.
$$
For a set $A$ of sequences of length $\ge n$, we put 
\begin{gather*}
A|_n=\{\aaa|_n:\aaa\in A\},
\qquad
A|^m=\{\aaa|^m:\aaa\in A\},
\intertext{and}
A|^m_n=\{\aaa|^m_n:\aaa\in A\}.
\end{gather*}
If $A$ is a set of finite sequences, $\ccc$ is  a
finite sequence, $B$ is a set of sequences, and 
$\ddd$ is  a sequence, then 
\begin{gather*}
\ccc B=\{\ccc\bbb:\bbb\in B\},
\qquad
AB = \{\aaa\bbb:\aaa\in A,\ \bbb\in B\},\\
\intertext{and}
A\ddd=\{\aaa\ddd: \aaa\in A\}.
\end{gather*}

Let $L$ be the set of all limit ordinals smaller 
than $\omega_1$, and let $S= \omega_1\setminus L$. 
We have $\omega_1 = \omega\cup \bigcup_{k\in \omega} (L+k)$, where $L+k = 
\{\alpha+k:  \alpha\in L\}$. 

Kulesza's space $Z\subset \omega_1^\omega$ is defined as 
\begin{align*}
Z= \{\aaa=a_0a_1\dots\in \omega_1^\omega:\ 
&\text{$a_0\in \omega_1\setminus L$, 
$a_k \in L$ for at most one $k\in \omega$,}\\ &\text{and if 
$a_k \in L$, then $a_{k+1}=a_k + k$  
and $a_{k+i}\in 
L+k$ for all $i\ge 2$}\}.  
\end{align*}

Kulesza proved that the space $Z$ with the topology induced by 
the topological product $\omega_1^\omega$ of countably many 
copies of the space $\omega_1$ with the usual order topology 
is metrizable and $\operatorname{{Ind}}Z=\dim Z=1$ (while, obviously, 
$\ind Z=0$)~\cite{Kul}. 

Kulesza did not give an explicit formula for a metric on $Z$, 
but he described base neighborhoods of the points of $Z$. They 
look as follows. 

For each limit ordinal $\alpha\in \omega_1$, we  fix 
an increasing sequence 
$\tilde\alpha_0\tilde\alpha_1\dots$ in $\omega_1$ with 
limit $\alpha$ and put $M_n(\alpha) = (\tilde \alpha_n, 
\alpha]$. 

Let $m\in \omega$. If a sequence $\aaa \in Z$ is such that 
$\aaa|_m\in S^m$, then we set 
$$
N_m(\aaa) = \{\bbb \in Z: \bbb|_m = 
\aaa|_m\}.  
$$
If $1\le k <m$ and $a_k\in L$, then 
$$
N_m(\aaa) = \{\bbb \in Z: \bbb|_k = 
\aaa|_k,\ \bbb|^{k+1}_{m}=\aaa|^{k+1}_{m},\
b_k\in 
M_m(a_k)\}.  
$$

The sets $N_m(\aaa)$ form 
a neighborhood base at the point $\aaa$ in the space $Z$. 

To prove the inequality $\dim Z>0$, Kulesza used 
the notion of full sets introduced by Fleissner in~\cite{Fle}. 

\begin{defn}[\cite{Fle}]
A set $T\subset \omega_1^{n}$ is said to be \emph{full} if 
$\{b_j: \bbb\in T,\ \bbb|_j=\aaa|_j\}$ is uncountable 
for any $\aaa \in T$ and $j<n$ (in particular, $T|_1$ 
is uncountable).  
\end{defn}

We say that a set $T\subset \omega_1^\omega$ is 
\emph{full} if $T|_n$ is full for all $n\in \omega$.

We need the following two combinatorial properties of full sets. 

\begin{lem}[{\cite[Lemma~6.4(b)]{Fle}}]\label{cl1}
If a set $T\subset \omega_1^n$ is full and $h\colon T\to 
\omega$, then $T$ contains a full subset on which $h$ 
is constant.  
\end{lem}

\begin{lem}\label{cl2}
If a set $T\subset \omega_1^\omega$ is full and 
$\{C_m:  m\in \omega\}$ is a family of sets such that $C_m\subset 
T|_m$ for $m\in \omega$ and, for any $\aaa \in T$, there exists 
an $n\in \omega$ for which $\aaa|_n\in C_n$, then  $C_t$ 
contains a full set {\upshape(}a subset of $T|_t\subset 
\omega_1^t${\upshape)} for some $t\in \omega$.  
\end{lem}

\begin{proof}
This lemma is similar to Lemma~6.4(a) from~\cite{Fle}. 
In~\cite{Fle}, the role of $T$ is played by $\omega_1^\omega$. 
There exists a natural bijection 
$$
\psi\colon [\omega_1]^{\le \omega}\to \bigcup_{n\le \omega}T|_n.
$$
It is constructed as follows. 
For all $n\in \omega$ and $\xxx\in T|_n$, we fix 
bijections $\varphi_{\xxx}\colon \omega_1\to \{y: \xxx y\in 
T|_{n+1}\}$ and put 
$$
\psi(\alpha_0\alpha_1\dots)=
\varphi(\alpha_0)\varphi_{\varphi(\alpha_0)}(\alpha_1)
\varphi_{\varphi(\alpha_0)\varphi_{\varphi(\alpha_0)}(\alpha_1)}(\alpha_2)\dots
$$
for any (finite or infinite) 
sequence 
$\alpha_0\alpha_1\dots \in 
[\omega_1]^{\le\omega}$. The map $\psi$ respects restrictions in the 
sense that if $\boldsymbol\alpha, \boldsymbol\beta\in 
[\omega_1]^{\ge n}$ and $\boldsymbol\alpha|_n=\boldsymbol\beta|_n$, 
then $\psi(\boldsymbol\alpha)|_n = \psi(\boldsymbol\beta)|_n$; 
moreover, $\psi (\omega_1^{n})= T|_n$.
The family $\{\psi^{-1}(C_m):m\in \omega\}$ has 
the properties 
$$
\psi^{-1}(C_m)\subset \omega_1^\omega|_m \quad
\text{for all $m\in \omega$}
$$
and 
$$
\text{for any $\boldsymbol \alpha\in \omega_1^\omega$, there exists 
an $n\in \omega$ such that  $\boldsymbol \alpha|_n\in 
\psi^{-1}(C_n)$.}
$$
According to~\cite[Lemma~6.4(a)]{Fle}, there exists a $t\in \omega$ 
for which $\psi^{-1}(C_t)$ contains a full set. For this $t$, 
$C_t$ contains a full set.  
\end{proof}

Levin~\cite{Lev} suggested a simple short proof 
of the inequality $\dim Z>0$ based on the notion of regular sets. 
We need the following modification of this notion. 

\begin{defn}
Let $U\subset Z\times Z$ be any set containing 
the diagonal. We say that a pair of sequences 
$(\xxx, \yyy)\in S^n\times S^n$ 
is \emph{$U$-regular} (or simply \emph{regular}, when it is clear 
what set $U$ is meant) if there exists a map 
(\emph{regulator}) $f\colon ([S]^{<\omega})^2\to \omega_1$ such that 
$(\xxx\aaa, \yyy\bbb)\in U$ 
whenever the sequences $\aaa, \bbb\in S^\omega$ satisfy the 
condition $a_i, b_i>f(\aaa|_i, \bbb|_i)$ for all $i\in\omega$ (in 
particular, $a_0, b_0> f(\varnothing)$).  
\end{defn}

Let $U$ be an arbitrary subset of $Z\times Z$ containing 
the diagonal. For $\aaa\in Z$, we put 
$$%
U(\aaa)= \{\bbb\in Z: (\aaa, \bbb)\in U\}.
$$%
The set $U^2$ is defined standardly as 
$$
U^2=\{(\aaa, \bbb): \text{there exists a $\ccc\in Z$ such that 
$(\aaa, \ccc)\in U$ and $(\ccc, \bbb)\in U$}\}.
$$
Thus, 
$$
U^2(\aaa)=\{\bbb\in Z: \text{there exists a $\ccc\in Z$ such 
that $(\aaa, \ccc)\in U$ and $(\ccc, \bbb)\in U$}\}.
$$

Suppose that $\{U_n:n\in \omega\}$ is a countable base 
for a uniformity on $Z$ generating the topology of the space $Z$. 
For each $\aaa\in S^\omega\subset Z$, fix $m_\aaa\ge 2$ 
for which $U^2_{m_\aaa}(\aaa)\subset N_2(\aaa)$. For $k\in 
\omega$, we set 
$$
C_k=\{\aaa|_k: \aaa\in S^\omega,\ m_\aaa\le k,\
N_2(\aaa)\supset U^2_{m_\aaa}(\aaa)\supset U_{m_\aaa}(\aaa)\supset 
N_k(\aaa)\}.
$$
Clearly, for any sequence $\aaa \in Z$, 
there exists a $k\ge m_\aaa$ for which $N_k(\aaa)\subset 
U_{m_\aaa}(\aaa)$ (because the sets $ U_{m_\aaa}(\aaa) $ are open 
and the $N_k(\aaa)$ form a base for the topology of $Z$ at the point 
$\aaa$).   Hence, for any sequence $\aaa \in S^\omega$, 
there exists a $k$ for which $\aaa|_k\in C_k$. By Lemma~\ref{cl2}, 
there exists a $t$ such that $C_t$ contains a full set 
(clearly, $t\ge 2$, because the sets $C_k$ are empty for $k<2$).  
Using Lemma~\ref{cl1}, we choose a number $m\in \omega$ and a 
full set $T\subset C_t$ such that $\min\{m_\aaa:\aaa|_t=a_0\dots 
a_{t-1}\}=m$ for any $a_0\dots a_{t-1}\in T$; note that 
$m\le t$ by the definition of $C_t$. We put $U=U_m$.  Our 
purpose is to show that $U\ne \overline U$.  
Suppose that $U= \overline U$.

\begin{rem}\label{rem30}
For any $\xxx\in Z$ 
such that $\xxx|_t\in T$, we have $U(\xxx)\subset N_2(\xxx)$.  
Indeed, by the definition 
of $T$, there exists an $\aaa\in S^\omega$ for which 
$\aaa|_t=\xxx|_t$, $N_2(\aaa)\supset 
U^2_{m_\aaa}(\aaa)\supset U_{m_\aaa}(\aaa)\supset N_t(\aaa)$, and 
$m_\aaa=m\le t$ (i.e., $U_{m_\aaa}=U$). Since $\xxx|_t = 
\aaa|_t\in S^t$, we have $\xxx\in N_t(\aaa)$. Therefore, $\xxx\in U(\aaa)$, and 
$U(\xxx)\subset U^2(\aaa)\subset N_2(\aaa)$. Since $t\ge m\ge 
2$ and $\xxx|_t = \aaa|_t(\in S^t)$,  it follows that 
$N_2(\aaa)=N_2(\xxx)$; thus, $U(\xxx)\subset N_2(\xxx)$.  
\end{rem}

\begin{rem}\label{rem31}
The pair $(x, x)$ is not $U$-regular for any $x\in C_t|_1$. 
Indeed, suppose that $x\in C_t|_1$, the pair $(x, x)$ is 
regular, and 
$f\colon ([S]^{<\omega})^2\to \omega_1$ is the 
corresponding regulator.  Since the set $C_t$ is full, 
we can find $a_1, a_2, \dots\,, b_1, b_2, \dots \in S$ such 
that 
\begin{gather*}
a_1\ne b_1, \qquad 
xa_1a_2\dots a_{t-1}, xb_1b_2\dots b_{t-1}\in C_t,\\ 
a_{1}, b_{1}> f(\varnothing), \qquad \text{and}\qquad
a_{i+1}, b_{i+1}>f(a_1\dots a_i, b_1\dots b_i)\quad\text{for 
all $i\ge 1$.}
\end{gather*}
Let $a_0=b_0=x$. We have $\aaa|_t, \bbb|_t\in C_t$. According 
to Remark~\ref{rem30}, $U(\aaa)\subset N_2(\aaa)$. However, by 
the definition of a regular pair, we also have $(\aaa, \bbb)\in U$, 
i.e., $\bbb\in U(\aaa)$.  Therefore, $\bbb \in N_2(\aaa)$, which 
is false, because $b_1\ne a_1$.  
\end{rem}

\begin{rem}\label{rem32}
On the other hand, for any pair $(\xxx, \yyy)\in U\cap (S^\omega\times 
S^\omega$) 
(in particular, for any pair $(\xxx, \xxx)$, where $\xxx \in S^\omega$), 
there exists an $n\in \omega$ such that the pair $(\xxx|_n, \yyy|_n)$ 
is regular. Indeed, since $U$ is open and the sets $N_k(\xxx)$ and 
$N_k(\yyy)$ form bases of neighborhoods of the points $\xxx$ and $\yyy$, 
it follows that there exists an $n\in \omega$ for which 
$$%
N_n(\xxx)\times N_n(\yyy)\subset U;
$$%
this means that 
$(x_0x_1\dots x_{n-1}\aaa, y_0y_1\dots y_{n-1}\bbb)\in U$ for any 
$\aaa$ and $\bbb$ from $Z$, not only for those satisfying 
the condition from the definition of regular pairs.  
\end{rem}

\begin{lem}\label{lemma31}
Suppose that $k>0$\textup; 
$\xxx = x_0\dots x_{k-1}, \yyy = y_0\dots y_{k-1} \in 
S^k$\textup; the pairs $(\xxx|_n, \yyy|_n)$ with $n\le k$ 
are not regular\textup; 
and there exists an uncountable set $S'\subset S$ 
such that the pair $(\xxx z, \yyy z)$ is regular for any $z\in S'$. 
Then there exists a number $l>0$, points $x_k, \dots, x_{k+l-1}, 
y_k, \dots, y_{k+l-1} \in S$, and an uncountable set $S''\subset S$ 
such that the pairs $(x_0\dots x_n, y_0\dots y_n)$ with $n<k+l$ 
are not regular and the pair  $(x_0\dots x_{k+l-1} z, 
y_0\dots y_{k+l-1} z)$ is regular for any $z\in S''$.  
\end{lem}

\begin{proof}
Let $C\subset L$ be an arbitrary closed unbounded 
set of limit ordinals. Take $c_0\in C$ and $z_0\in S'$ for 
which $z_0> c_0$. By assumption, the pair $(\xxx z_0, \yyy z_0)$ 
is regular; let $f_0$ be the corresponding regulator. Take 
$c_1\in C$ such that $c_1>\max\{f_0(\varnothing), z_0\}$ and 
$z_1\in S^1$ such that $z_1> c_1$. By assumption, the pair $(\xxx z_1, 
\yyy z_1)$ is regular; let $f_1$ be the corresponding regulator. 
Suppose that we made $n$ steps, i.e., chose ordinals $c_{n-1}\in 
C$ and $z_{n-1}\in S'$ and a regulator $f_{n-1}$.  At the $(n+1)$th 
step, we take $c_n\in C$ and $z_n\in S'$ such that 
$$
c_n> \max\{f_{n-1}(\varnothing), z_{n-1}\}
\qquad \text{and}\qquad 
z_n> c_n,
$$
and choose a map $f_n$ witnessing the regularity 
of the pair $(\xxx z_n, \yyy z_n)$.

As a result, we obtain an increasing sequence 
of elements  of $C$. Let $c=\sup \{c_n:n\in \omega\}$. We have 
$c\in C$, because $C$ is closed. Moreover, for any $n\in \omega$, 
the pair $(\xxx z_n, \yyy z_n)$ is regular, $f_n$ is the 
corresponding regulator, and $c+k> c>f_n(\varnothing)$.  
Therefore, if $\aaa\in S^\omega$ is a sequence such that 
\begin{equation} \label{eq31}
a_i>\sup \{f_n((c{+}k)\aaa|_i, (c{+}k)\aaa|_i):n\in \omega\} 
\quad\text{for all $i\in \omega$}, 
\end{equation} 
then $(\xxx z_n(c{+}k)\aaa, \yyy z_n(c{+}k)\aaa)\in U$. 

Recall that $c=\sup \{c_n:n\in \omega\}=\sup\{z_n:n\in 
\omega\}$; thus, any neighborhood in $Z\times Z$ of any point of the form 
$(\xxx c (c{+}k)\aaa, \yyy c (c{+}k)\aaa)$ contains the point 
$(\xxx z_n (c{+}k)\aaa, \yyy z_n (c{+}k)\aaa)$ for some $n$. 
Therefore, if a sequence $\aaa$ satisfies condition~\eqref{eq31}, 
then
$$
(\xxx c (c{+}k)\aaa, \yyy c (c{+}k)\aaa)\in \overline U=U. 
$$

Clearly, the set of sequences $\aaa\in S^\omega$ 
satisfying~\eqref{eq31} is full. 

Thus, any closed unbounded set of limit ordinals 
contains a point $c\in L$ for which there exists a full 
set $Y_c\subset S^\omega$ such that 
$$
(\xxx c (c{+}k)\zzz, \yyy c (c{+}k)\zzz)\in U
\quad\text{for any $\zzz\in Y_c$}.
$$
Therefore, the set $L'$ of such points $c$ is stationary. 

Since $U$ open, it follows that, for any $c\in L'$ and 
$\zzz\in Y_c$, there exists an $n=n(\zzz, c)> k+2$ such that 
$$
N_n(\xxx c(c{+}k)\zzz)\times N_n(\yyy c(c{+}k)\zzz)\subset U.
$$
For $m\in \omega$ and $c\in L'$, we set 
$$
C_m(c)=\{\zzz\in Y_c: n(\zzz, c)=m\}|_m.
$$
For any $c\in L'$, using Lemma~\ref{cl2} and the definition of the 
neighborhoods of the form $N_n(\aaa)$, we can find an  
$m_c>0$ and a full set $Y'_c\subset Y_c|_{m_c-k-2}$ such 
that 
\begin{equation}\label{eq32}
(\xxx \mu (c{+}k)\zzz\aaa, \yyy \nu (c{+}k)\zzz\bbb)\in U\qquad
\text{for any $\mu, \nu\in M_{m_c}$,  $\zzz\in Y'_c$, and 
$\aaa, \bbb\in Z|^{m_{c}}$}.  
\end{equation}
Using the pressing down lemma, we choose 
a stationary subset $L''$ of the stationary set $L'$ such that 
$$
\tilde c_{m_c}=\beta \quad\text{for all $c\in L''$},
$$
where $\beta$ is a countable ordinal (here the $\tilde c_n$ 
are the ordinals converging to $c$ that are used in the definition 
of the sets $M_n(c)$ involved in the definition of the neighborhoods 
$N_n(\xxx c(c{+}k)\zzz)$ and $N_n(\yyy c(c{+}k)\zzz)$). 

Suppose that the pairs $(\xxx x, \yyy y)$ are regular for any 
$x,y> \beta$ from $S$ and $f_{xy}$ are the corresponding 
regulators.  Then the pair $(\xxx, \yyy)$ itself is regular: 
the corresponding regulator is defined by 
$$
f(\varnothing)=\beta, 
\qquad
f(x\aaa, y\bbb)=
\begin{cases}
f_{xy}(\aaa, \bbb)& \text{if $x,y> \beta$},\\
0& \text{if $x\le \beta$ or $y\le \beta$}.
\end{cases}
$$

The pair  $(\xxx, \yyy)$ is not regular by assumption; 
hence there exist $x_k, y_k> \beta$ for which the pair $(\xxx 
x_k, \yyy y_k)$ is not regular. 

If the pairs $(\xxx x_k (c{+}k), \yyy y_k (c{+}k))$ are regular for all 
$c\in L''$ such that $c> x_k, y_k$, then we can 
set $l=1$ and $S''=L''+k$. Otherwise, we take $c\in 
L''$ for which $c>x_k$, $y_k>\beta$, and the pair 
$(\xxx x_k (c{+}k), \yyy y_k (c{+}k))$ 
is not regular. Condition~\eqref{eq32} implies 
\begin{equation}\label{eq33}
(\xxx x_k (c{+}k)\zzz\aaa, \yyy y_k(c{+}k)\zzz\bbb)\in U
\quad\text{for any $\zzz\in Y'_c$ and $\aaa, \bbb\in Z|^{m_{c}}$}.
\end{equation}

The set $Y'_c$ is full; hence $Y'_c|_1$ is uncountable. 
If the pairs 
$(\xxx x_k (c{+}k)z,\allowbreak \yyy y_k (c{+}k)z)$ are regular 
for all $z\in Y'_c|_1$, then we have obtained what is required. 
Otherwise, we take 
$z_0\in Y'_c|_1$ for which the pair $(\xxx x_k (c{+}k)z_0, \yyy y_k 
(c{+}k)z_0)$ is not regular.  Relation~\eqref{eq32} implies 
$$
(\xxx x_k (c{+}k)z_0\zzz\aaa, \yyy y_k(c{+}k)z_0\zzz\bbb)\in U\qquad
\text{for any $\zzz$ such that $z_0\zzz\in Y'_c$ and 
any $\aaa, \bbb\in Z|^{m_{c}}$}.
$$
The set $Y'_c$ is full and $z_0\in Y'_c|_1$; hence the set $\{z: 
z_0z\in Y'_c|_2\}$ is uncountable.  If the pairs $(\xxx x_k 
(c{+}k)z_0z,\yyy y_k (c{+}k)z_0z)$ are regular for all 
$z\in S$ such that $z_0z\in Y'_c|_2$, then we have obtained what is 
required.  Otherwise, we continue the construction. Sooner or later, 
the procedure will terminate:  we shall find either an 
$n<m_c-k-4$ such that the pairs $(\xxx x_k (c{+}k)z_0\dots 
z_{n-1}z,\allowbreak 
\yyy y_k (c{+}k)z_0\dots z_{n-1}z)$ are regular for all 
$z$ with $z_0\dots z_{n-1}z\in Y'_c|_{n+1}$ or 
$z_0\dots z_{m_c-k-4}\in Y'_c|_{m_c-k-3}$ such that all pairs $(\xxx 
x_k (c{+}k)z_0\dots z_{m_c-k-4}|_n, \yyy y_k (c{+}k)z_0\dots 
z_{m_c-k-4}|_n)$, where $n\le m_c$, are not regular. In the latter 
case, the pair 
$$
(\xxx x_k (c{+}k)z_0\dots 
z_{m_c-k-4}z, \yyy y_k (c{+}k)z_0\dots z_{m_c-k-4}z)
$$
is regular for 
any $z$ such that $z_0\dots z_{m_c-k-4}z\in Y'_c$ (and there are 
uncountably many such $z$, because $Y'_c$ is full) by 
virtue of~\eqref{eq33}.  
\end{proof}

Take any point $x_0\in C_t|_1$ (the set $C_t$ was 
defined before Remark~\ref{rem30}). According to 
Remark~\ref{rem30}, the pair $(x_0, x_0)$ is not regular. If there 
exists an $x\in S$ for which the pair $(x_0x, x_0x)$ is not regular, 
then we take this $x$ for $x_1$. Suppose that we have constructed 
a sequence $x_0x_1\dots x_{n-1}\in S^n$ such that the pairs 
$(x_0x_1\dots x_{i-1}, x_0x_1\dots x_{i-1})$ are not regular for any 
$i\le n$.  If there exists an $x\in S$ for which the pair 
$(x_0x_1\dots x_{n-1}x, x_0x_1\dots x_{n-1}x)$ is not regular, then 
we take this $x$ for $x_n$.  The construction cannot be continued 
infinitely long (otherwise, we 
shall obtain a sequence $\xxx\in S^\omega$ such that the pair 
$(\xxx|_n, \xxx|_n)$ is not regular for any $n\in 
\omega$, whose existence contradicts Remark~\ref{rem32}). Thus, 
sooner or later, we shall obtain a sequence $x_0\dots 
x_{k-1}\in S^k$ such that the pair $(x_0x_1\dots x_{i-1}, x_0x_1\dots 
x_{i-1})$ is not regular for any $i\le k$ but all 
pairs $(x_0x_1\dots x_{k-1}x, x_0x_1\dots x_{k-1}x)$, where $x\in S$, 
are regular. 

We set $y_0\dots y_{k-1}=x_0\dots x_{k-1}$. Applying 
Lemma~\ref{lemma31} to the pair $(x_0\dots x_{k-1},\allowbreak 
y_0\dots y_{k-1})$, we obtain a pair  $(x_0\dots x_{k'-1}, y_0\dots 
y_{k'-1})$ such that $k'>k$, the pair $(x_0x_1\dots x_{i-1}, 
y_0y_1\dots y_{i-1})$ is not regular for any $i\le k'$, but all pairs 
$(x_0x_1\dots x_{k'-1}z, y_0y_1\dots y_{k'-1}z)$, where $z$ belong to 
some uncountable set $S'\subset S$, are regular. Repeatedly applying 
Lemma~\ref{lemma31}, we shall extend the sequences in this pair. In 
the end, we shall obtain a sequence $\xxx,\yyy\in S^\omega$ such that, 
for any $n>0$,  
the pair $(\xxx|_n, \yyy|_n)$ is not regular but  
there exists an $m\ge n$ and an uncountable set 
$S_m\subset S$ such that all pairs $(\xxx|_m z, \yyy|_m z)$, where 
$z\in S_m$, are regular. 

Take any $n\in \omega$ and consider the neighborhood $N_n(\xxx)\times 
N_n(\yyy)$ of the pair $(\xxx, \yyy)$ in $Z\times Z$. Suppose that $m\ge n$ 
and $z\in S$ are such that the pair $(\xxx|_m z, \yyy|_m z)$ 
is regular.  This means that  
$(\xxx|_m z\aaa, \yyy|_m z\bbb)\in U$ 
for some $\aaa, \bbb\in Z|^{m+1}$. Clearly, $(\xxx|_m z\aaa, 
\yyy|_m z\bbb)\in N_n(\xxx)\times N_n(\yyy)$. Thus, $N_n(\xxx)\times 
N_n(\yyy)\cap U\ne \varnothing$ for any $n\in \omega$, and hence 
$(\xxx, \yyy)\in \overline U=U$. Remark~\ref{rem32} 
implies that the pair $(\xxx|_n, \yyy|_n)$ must be regular for some $n$. 
This contradiction shows that $\overline U\ne U$. 

\section
{Concluding Remarks}

We have considered two metrizable spaces with noncoinciding 
dimensions, Mrowka's and\break 
Kulesza's, and 
shown that one of them can be embedded in 
a zero-dimensional metrizable group  and the other cannot. The 
natural question arises: What properties of Kulesza's space obstruct 
its embedding into a zero-dimensional metrizable group? The 
most manifest difference between Mrowka's and Kulesza's spaces 
is that Kulesza's space 
is metrizable by a complete metric. This suggests the 
conjecture that a space metrizable by a complete metric can be 
embedded in a zero-dimensional metrizable group only if it 
is strongly zero-dimensional.  
This conjecture is based not only on purely formal grounds
but also on some intuitive reasons; in this author's opinion, it is 
fairly likely. Even more likely is the following auxiliary 
conjecture: If $(X, \rho)$ is a metric space with complete metric 
$\rho$, $A_\rho(X)$ is the free group of $X$ metrized by the 
Graev extension of $\rho$, and $\ind A_\rho(X)=0$, 
then $\dim X=0$. 

It is also unclear how the dimension 
of metrizable groups behaves under completion\footnote{This 
question is difficult even for general 
topological spaces.  Thus, Mrowka's space $\nu\mu_0$ has 
a zero-dimensional completion under the 
continuum hypothesis~\cite{Mro1}; however, Mrowka also proved 
that the assertion that the small inductive dimension of all metric 
completions of $\nu\mu_0$ is larger than zero is possibly 
consistent~\cite{Mro1}, i.e., it holds 
under a certain set-theoretic assumption 
whose consistency with ZFC is very likely.}.  
It is only clear that the free and free Abelian 
groups with Graev metrics (as well 
as the metrizable groups of the form $(A(X), 
\mathcal T_{\boldsymbol{\mathcal U}})$ described in the first section,  
into which we can embed zero-dimensional 
metrizable spaces) are never complete; we can always 
construct a fundamental sequence consisting of words with unboundedly 
increasing lengths, which converges to no word of finite 
length.\footnote{More details on topologies on free 
groups (including 
the Graev metric topology) can be found in~\cite{Sipasurvey}.} 

We conclude this paper with several questions. 

\begin{prob}
Is it true that if the uniformity  generated by a metric $\rho$ 
on a set $X$ has a countable base consisting of open-and-closed 
sets, then the free (Abelian) group of $X$ metrized with the Graev 
extension of $\rho$ is zero-dimensional?
\end{prob}

\begin{prob}
Does there exist a complete metric group with noncoinciding 
dimensions $\ind$ and $\dim$?  
\end{prob}

\begin{prob}
Is it true that any complete metric space 
which can be embedded into a  zero-dimensional metrizable group 
is strongly zero-dimensional?  
\end{prob}

\begin{prob}
Is it true that 
if $(X, \rho)$ is a complete metric space with 
metric $\rho$, $G_\rho(X)$ is the free (Abelian) 
group of $X$ metrized by 
the Graev extension of $\rho$, and $\ind 
G_\rho(X)=0$, then $\dim X=0$?  
\end{prob}

\begin{prob}
Is it true that 
if $(X, \rho)$ is a metric space with 
metric $\rho$, $G_\rho(X)$ is the free (Abelian) 
group of $X$ metrized by 
the Graev extension of $\rho$, and the completion 
of $G_\rho(X)=0$ is zero-dimensional, then 
$\dim X=0$?  What if the metric $\rho$ is complete?  
\end{prob}

\begin{prob}
How large can the gap between the dimensions $\ind$ and 
$\dim$ of a metrizable group be? 
What values can the dimension $\dim$ of a metrizable 
topological group $G$ with $\ind G=0$ take?  
\end{prob}

\begin{prob}
Let $(\nu\mu_0(A), \rho)$ be 
Mrowka's space described in the second section with 
a metric $\rho$ generating 
the uniformity with a clopen base described in the same section, and 
let be $G$ the metrizable group with $\ind G= 0$ into which 
$\nu\mu_0(A)$ is embedded by Theorem~\ref{theorem1.1}. 
\begin{enumerate}
\item[(a)]
Find $\dim G$; 
\item[(b)]
Find 
$\ind G_\rho(\nu\mu_0(A))$ and $\dim G_\rho(\nu\mu_0(A))$, 
where $G_\rho(\nu\mu_0(A))$ is the free (Abelian) group of $\nu\mu_0(A)$ 
metrized by the Graev extension of the metric $\rho$. 
\end{enumerate}
\end{prob}


\end{document}